\documentclass[11pt]{article}
\usepackage{paper}
\usepackage[black,external]{causets}
\tikzcausetsset{
  text scale=1,
  text font=\scriptsize,
  link width=0.5pt,
  wild/.style={draw, fill=#1},
  wild/.default={white},
  left label/.style={label=#1, replaced labels={left=0.25ex, at=west}}
}
\colorlet{sel1}{cyan}
\colorlet{sel2}{magenta}
\colorlet{sel3}{yellow!75!black}
\graphicspath{{figures/}} 

\newcommand*{\npcauset}[2][]{\pcauset[name=#1]{#2}}
\newcommand*{\nrcauset}[3][]{\rcauset[name=#1]{#2}{#3}}
\newcommand*{\ncauset}[3][]{\causet[name=#1]{#2}{#3}}
\newcommand*{\smallTwoChain}{\pcauset[small, name=2Chain_small]{1,2}}
\newcommand*{\smallSingleton}{\pcauset[small, name=Singleton_small]{1}}

\title{Local symmetries in partially ordered sets}
\author{Christoph Minz\footnote{Institut f\"ur theoretische Physik, Universit\"at Leipzig, Br\"uderstra\ss e 16, 04103 Leipzig, Germany; e-mail: \url{christoph.minz@itp.uni-leipzig.de}}}
\date{June 17, 2024}

\usepackage{csquotes}  
\usepackage[%
backend=biber,
style=alphabetic,  
citestyle=alphabetic,  
sorting=nyt,
maxnames=10,  
giveninits=true,  
eprint=true,  
doi=false,  
url=false  
]{biblatex}
\begin{document}

\maketitle

\begin{abstract}
  Partially ordered sets (posets) play a universal role as an abstract structure in many areas of mathematics.
  For finite posets, an explicit enumeration of distinct partial orders on a set of unlabelled elements is known only up to a cardinality of 16 (listed as sequence A000112 in the OEIS), but closed expressions are unknown.
  
  By considering the automorphisms of (finite) posets, I introduce a formulation of local symmetries.
  These symmetries give rise to a division operation on the set of posets and lead to the construction of symmetry classes that are easier to characterise and enumerate.
  Furthermore, we consider polynomial expressions that count certain subsets of posets with a large number of layers (a large height).
  As an application in physics, local symmetries or rather their absence helps to distinguish causal sets (locally finite posets) that serve as discrete spacetime models from generic causal sets.
\end{abstract}

\section{Introduction}
\label{sec:Introduction}

Counting is perhaps the most basic mathematical concept and we may find unexpected connection between different mathematical structures whenever they have an identical enumeration.
A fairly ubiquitous structure in various fields of mathematics are partial orders. 
For example, all subspaces of a von Neumann algebra are partially ordered under inclusion, which implies a partially ordered set of the corresponding subspace projections.
Another example are the incidence relations between edges and vertices of a graph, or even the incidences of any (abstract) geometric shape with edges, vertices and higher faces has a representation as a partially ordered set.

The number of distinct partially ordered sets (posets) on a set of $n$ unlabeled elements is listed as sequence A000112 in the OEIS~\cite{OEIS:A000112:2023} (identical to the number of $T_0$ topological spaces), but only the first sixteen terms of the sequence have been determined to date \cite{BrinkmannMcKay:2002}, and enumeration results for certain subsets have been studied~\cite{Klarner:1969,Butler:1972}.
Although the asymptotic behaviour of the sequence is known~\cite{KleitmanRothschild:1975,Dhar:1980,ProemelStegerTaraz:2001GlobalStructure,ProemelStegerTaraz:2001Evolution}, finding analytic expression for the number of partial orders remains challenging.
This work shall provide a classification of posets by analysing repeating substructures forming local symmetries and making use of the concept of retracts~\cite{DuffusRival:1979,DuffusRival:1981}.
Although some of the findings might already be known, I use the characterisation of posets with and without symmetries to formulate explicit expressions for partial enumerations, and I illustrate an application in the quantum gravity programme of causal set theory.

In physics, the causal structure of a spacetime manifold is a partial order, and as an approach to quantum gravity, causal set theory is based on discrete counterparts to spacetime manifolds where any interval is finite and it is hypothesised that continuum spacetimes are corresponding to or ``emerging'' from causal sets~\cite{BombelliLeeMeyerSorkin:1987,BrightwellGregory:1991,BrightwellHensonSurya:2007}.
Any local subset (interval) of a causal set is a finite poset.
While causal sets that correspond to a given spacetime manifold are usually generated via a random process (called sprinkling and mathematically reviewed in~\cite{FewsterHawkinsMinzRejzner:2021}), the central idea of the causal set programme is to describe spacetime manifolds as continuous counterparts of the fundamentally discrete structures.
The asymptotic description of generic posets in~\cite{KleitmanRothschild:1975}, however, showed that generic causal sets are most likely of a type that is very different from the causal sets obtained by a sprinkling process.
Recently, it was shown that the typical Kleitman-Rothschild orders could be sufficiently suppressed in the ``sum-over-history'' approach leaving us with those causal sets that more likely represent discrete spacetimes~\cite{CarlipCarlipSurya:2023}.
Below, we consider properties to help distinguish between generic posets and those that are more relevant for causal set theory.

The work is structured as follows.
First, we recall definitions for the characterization of elements, subsets and entire posets.
The definition of local symmetries in posets leads us to the classification of posets into symmetry classes, which often have explicit enumerations.
I give closed expressions for some of these classes, and also for the enumeration of posets with a large number of layers.
Finally, I comment on the application in causal set theory.

\section{Notation and terminology}
\label{sec:Notation}

In this section, we recall some basic notation and terminology that is being used throughout this work. 
Some readers might be familiar with partial orders, though, I recommend to skim over this section, at least, as I also include some new or less common terms here. 

\begin{definition}
\label{def:Poset}
  A \defof{partially ordered set} (\emph{poset}) is a set $P$ with a partial ordering $\leq$ that is (for all $a, b, c \in P$)
  \begin{subequations}
  \begin{align}
      \text{transitive:}&
  &
      \left( a \leq b \text{~and~} b \leq c \right)
      \limp a
    &\leq c
    \eqend{,}
  \nexteq
      \text{antisymmetric:}&
  &
      \left( a \leq b \text{~and~} b \leq a \right)
      \limp a
    &= b
    \eqend{,}
  \nexteq
      \text{reflexive:}&
  &
      a
    &\leq a
    \eqend{.}
  \end{align}
  \end{subequations}
\end{definition}
\begin{definition}
\label{def:OppositePoset}
  The \defof{opposite poset} of a poset $(P, \leq)$ is the poset $(P, \binop{\leq})$ such that
  \begin{align}
  \label{eq:OppositePoset}
        \forall a, b
    &\in P:
  &
        b \binop{\leq} a 
    &\liff a \leq b
    \eqend{.} 
  \end{align}
  A poset is \defof{self-opposite} if it is isomorphic to its opposite. 
\end{definition}
For short, I denote a poset by its set $P$ while the ordering $\leq$ is either kept implicit or can be inferred from the context.
Recall that any isomorphism between (sub)posets $A$ and $B$, written as $A \ito B$, is order-preserving and order-reflecting.
In this work, I consider partial orders on sets of unlabelled elements only, and the set of all posets (up to isomorphisms) is denoted by $\mathfrak{P}$.
The set of all isomorphisms from a poset $P \in \mathfrak{P}$ to itself forms the automorphism group $\Aut(P)$.
A subset $A$ of the poset $P$ is \emph{fixed} under an automorphism $\alpha \in \Aut(P)$ if $\forall a \in A: \alpha(a) = a$.
Let $\Sigma(\alpha) \subseteq P$ denote the set of all elements that are not fixed under $\alpha$.

Any subset of a poset is also a poset with the restricted partial order.
A \emph{chain} is a totally ordered poset (or total order) and an \emph{antichain} is a poset that contains only pairwise separated elements.
\begin{definition}[Fences]
\label{def:Fences}
  An \defof{$n$-fence} $F_n$ is a poset of $n$ elements $a_i$ that are partially ordered in alternation,
  \begin{subequations}
  \label{eq:Fences}
  \begin{align}
    a_1 < a_2 > a_3 < \dots a_n
  &\qquad\text{(lower fence) or}\qquad
  \nexteq
    a_1 > a_2 < a_3 > \dots a_n
    &\qquad\text{(upper fence)}
    \eqend{.}
  \end{align}
  \end{subequations}
  A fence with an even number of elements is both upper and lower.
\end{definition}

\begin{definition}[Element properties]
\label{def:ElementSubsets}
  Let $P$ be a poset. 
  The sets of elements \defof{preceding and succeeding $a \in P$} are denoted by 
  \begin{subequations}
  \begin{align}
  \label{eq:PrecedingElements}
      J^-(a)
    &:= \left\{
        c \in P
      \bypred
        c \leq a
      \right\}
    \eqend{,}
  \nexteq
  \label{eq:SucceedingElements}
      J^+(a)
    &:= \left\{
        c \in P
      \bypred
        a \leq c
      \right\}
    \eqend{,}
  \end{align}
  \end{subequations}
  respectively.
  For the strict versions, write $J^\pm_*(a) = J^\pm(a) \setminus \{a\}$.
  An element $a$ is \defof{covered by} $b$ and $b$ is \defof{covering} $a$ written as 
  \begin{align}
  \label{eq:LinkedElements}
      a \lessd b
    &:\liff \left( a < b
    \text{~and~} \nexists c \in P: a < c < b \right)
    \eqend{.}
  \end{align}
  The sets of all \defof{linked} elements are those covered by and covering $a$, denoted as $L^\pm(a) \subseteq J^\pm(a)$, respectively.
  Elements $a$ and $b$ are \defof{fence connected} if there exists a fence subset $F \subseteq P$ such that $a \in F$ and $b \in F$.
  The elements are \defof{separated} (also \emph{independent} or \emph{incomparable}) if $a \nleq b$ and $a \ngeq b$.
\end{definition}
Any single element as well as any two distinct but related elements in a poset are trivially fence connected.
For a poset $P$, a \emph{component} is a subset $A \subseteq P$ where all elements in $A$ are fence connected and no element of $P \setminus A$ is fence connected to any element in $A$.
\begin{remark}
  All separated element pairs of a poset form a symmetric binary relation that is an irreflexive locality relation~\cite{GuoPaychaZhang:2023}.
\end{remark}

\begin{definition}[Subset properties]
\label{def:PosetSubsets}
  A subset $A$ of a poset $P$ is \defof{connected} if all pairs of elements $a, b \in A$ are fence connected.
  It is also \defof{linked} if
  \begin{align}
      \forall a, b \in A:
      a < b
    &\limp \Bigl(
        \bigl( a \lessd b \bigr)
        \text{~or~}
        \bigl( \exists c \in J^+(a) \cap J^-(b): c \in A \bigr)
      \Bigr)
    \eqend{.}
  \end{align} 
  I call a subset $A$ \defof{maximally ordered in $P$} if the number of ordered pairs 
  \begin{align}
      \Bigl\lvert \left\{ (a, b) \in A \times A \bypred a < b \right\} \Bigr\rvert
  \end{align}
  is maximal among all subsets of $P$ with cardinality $\lvert A \rvert$. 
  Two proper subsets $A, B \subset P$ are \defof{separated} if all pairs of elements $a \in A$ and $b \in B$ are separated.
\end{definition}

A maximally ordered subset of a poset is a chain or it has the same \emph{height} as the poset in the following sense.
\begin{definition}
\label{def:Layers}
  The \defof{number of layers} of a poset $P$, written as $\ell(P)$, is the length (cardinality) of its longest chain.
  The \defof{layer of an element} $a \in P$ is the number of layers of its preceding subset, $\ell_{P}(a) := \ell\bigl( J^-(a) \bigr)$.
\end{definition}
This layer number function will be used in the classification of posets below.
Properties like the number of layers are often easier to study with a representation of a poset as a Hasse diagram, the (directed) graph with vertices for all elements and edges (links) for all covering relations (pointing towards the top of the page).
Such diagrams are plotted with the \LaTeX-package ``causets''~\cite{Minz:2020}, for example, all posets with a cardinality $n \in \{1, 2, 3, 4\}$ are
\begin{subequations}
\label{eq:FinitePosets}
\begin{align}
    \mathfrak{P}_1
  &= \left\{
      \npcauset[Singleton]{1}
    \right\}
  \eqend{,}
\nexteq
    \mathfrak{P}_2
  &= \left\{
      \left( \npcauset[2Antichain]{2,1} \right),
      \npcauset[2Chain]{1,2}
    \right\}
  \eqend{,}
\nexteq
    \mathfrak{P}_3
  &= \left\{
      \left( \npcauset[3Antichain]{3,2,1} \right),
      \left( \npcauset[2Chain_Singleton]{3,1,2} \right),
      \npcauset[Wedge]{2,1,3},
      \npcauset[Vee]{1,3,2},
      \npcauset[3Chain]{1,2,3}
    \right\}
  \eqend{,}
\nexteq
    \mathfrak{P}_4
  &\pickindent{= \Bigl\{}
      \left( \npcauset[4Antichain]{4,3,2,1} \right),
      \left( \npcauset[2Chain_2Antichain]{4,3,1,2} \right),
      \left( \npcauset[Wedge_Singleton]{4,2,1,3} \right),
      \left( \npcauset[Vee_Singleton]{4,1,3,2} \right),
      \left( \npcauset[2Chain_2Chain]{3,4,1,2} \right),
      \npcauset[3Wedge]{3,2,1,4},
      \npcauset[3Vee]{1,4,3,2},
      \npcauset[N]{2,4,1,3},
      \npcauset[2Gonal]{2,1,4,3},
  \eqbreakr
      \left( \npcauset[3Chain_Singleton]{4,1,2,3} \right),
      \npcauset[Downhook]{3,1,2,4},
      \npcauset[Uphook]{1,4,2,3},
      \npcauset[DownY]{2,1,3,4},
      \npcauset[UpY]{1,2,4,3},
      \npcauset[2Diamond]{1,3,2,4},
      \npcauset[4Chain]{1,2,3,4}
    \Bigr\}
  \eqend{.}
\end{align}
\end{subequations}
For better readability, I parenthesize posets that are constructed from components by parallel composition, which is an operation on the class of posets.
\begin{definition}
\label{def:Composition}
  The \defof{series/parallel composition} of two posets $P$ and $Q$, denoted as $P \vee Q$ and $P \sqcup Q$, respectively, is a poset $R$ that is the disjoint union of $P$ and $Q$ such that any two elements both from $P$ or both from $Q$ have the same ordering in $R$. 
  Any two elements $p \in P$ and $q \in Q$ fulfill $p < q$ in $R$ for the series composition, respective $a \nleq b$ and $b \nleq a$ for the parallel composition.
\end{definition}

For the application in physics (mostly in the final part), recall the following definition.
\begin{definition}
\label{def:CausalSet}
  A \defof{causal set} (\emph{causet}) is a poset $P$ with an order that is (for all $a, b \in P$)
  \begin{align}
  \label{eq:CausalSet.Finiteness}
      \text{locally finite:}&
    &
      \Bigl\lvert J^+(a) \cap J^-(b) \Bigr\rvert
    &< \infty
    \eqend{.}
  \end{align}
\end{definition}
Whenever the partial order describes causality of a physical system, the sets $J^\mp(a)$ are also referred to as \emph{past and future} of an element $a$, respectively.

\section{Local symmetries in posets}
\label{sec:LocalSymmetries}

As a pre-cursor to the definition of local symmetries, note the following fact about subsets of strictly preceding and succeeding elements. 
\begin{lemma}[Common links]
\label{lma:CommonLinks}
  Let $P$ be a poset. 
  For all elements $a, b \in P$, 
  \begin{align}
  \label{eq:CommonLinkedSubset}
      \Bigl( J^\pm_*(a) = J^\pm_*(b) \Bigr)
    &\liff \Bigl( L^\pm(a) = L^\pm(b) \Bigr)
    \eqend{.}
  \end{align}
\end{lemma}
\proof
If the left hand side holds, then we have, for every $c \in J^\pm_*(a) = J^\pm_*(b)$, 
\begin{align}
    \bigl( c \in L^\pm(a) \bigr)
  \liff \bigl(
      J^\pm_*(a) \cap J^\pm_*(c)
      = \emptyset
    \bigr)
  &\liff \bigl(
      J^\pm_*(b) \cap J^\pm_*(c)
      = \emptyset
    \bigr)
  \liff \bigl( c \in L^\pm(b) \bigr)
  \eqend{,}
\end{align}
so $L^\pm(a) = L^\pm(b)$. 
If the right hand side holds, then the left hand side follows by transitivity of the partial order. 
\qed

\subsection{Symmetric elements in antichains}
\label{subsec:LocalSymmetries.Element}

For the introduction of local symmetries, we start with the simplest case. 

\begin{definition}
\label{def:ElementSymmetry}
  Let $P$ be a poset. 
  Two elements $a, b \in P$ are \defof{singleton-symmetric} if 
  \begin{align}
  \label{eq:Symmetry.Element}
        L^\pm(a)
    &= L^\pm(b)
    \eqend{.}
  \end{align}
\end{definition}
For any pair of singleton-symmetric elements $a, b \in P$, also $J^\pm_*(a) = J^\pm_*(b)$ by \autoref{lma:CommonLinks}. 
Immediate from the definition is that the singleton-symmetry is an equivalence relation (reflexive, symmetric and transitive).
Since each pair of elements in an equivalence class has the same order relations to the elements of the other equivalence classes, the quotient inherit the given partial order of the original poset $P$. 
The quotient of $P$ by all singleton-symmetries is isomorphic to a \emph{retracted} poset or \emph{retract}, written as $P \oslash \npcauset[Singleton]{1}$. 

\begin{example}
\label{eg:Antichains.Singleton}
  All antichain posets with more than one element are singleton-symmetric, 
  \begin{align}
  \label{eq:Example.Antichains.Singleton}
      \npcauset[Singleton]{1}
    &= \left( \npcauset[2Antichain]{2,1} \right)
      \oslash \npcauset[Singleton]{1}
    = \left( \npcauset[3Antichain]{3,2,1} \right)
      \oslash \npcauset[Singleton]{1}
    = \left( \npcauset[4Antichain]{4,3,2,1} \right)
      \oslash \npcauset[Singleton]{1}
    = \dots
    \eqend{.}
  \end{align}
  Note that the subgroup of automorphisms of an $n$-antichain is isomorphic to the permutation group $S_n$. 
\end{example}
\begin{example}
\label{eg:CompleteBipartite.Singleton}
  Similarly, the 2-chain is the retracted poset of 
  \begin{align}
  \label{eq:Example.CompleteBipartite.Element}
        \npcauset[2Chain]{1,2}
    &= \npcauset[Wedge]{2,1,3}
        \oslash \npcauset[Singleton]{1}
    = \npcauset[Vee]{1,3,2}
        \oslash \npcauset[Singleton]{1}
    = \npcauset[3Wedge]{3,2,1,4}
        \oslash \npcauset[Singleton]{1}
    = \npcauset[2Gonal]{2,1,4,3}
        \oslash \npcauset[Singleton]{1}
    = \npcauset[3Vee]{1,4,3,2}
        \oslash \npcauset[Singleton]{1}
    = \dots
    \eqend{.}
  \end{align}
\end{example}
\begin{example}
\label{eg:ParallelComposition.Singleton}
  The quotient by the singleton-symmetry distributes over parallel compositions of non-singleton posets (and similarly for series compositions) as in 
  \begin{subequations}
  \begin{align}
  \label{eq:Example.ParallelComposition}
        \left(
          \nrcauset[3CrownBrokenLow]{5,4,2,7,1,6,3}{2/6}
        \sqcup \npcauset[2CrownLinked3Wedge]{7,6,5,2,8,1,4,3}
        \sqcup \npcauset[2Chain3Chain]{2,1,5,4,3}
        \right)
        \oslash \npcauset[Singleton]{1}
    &= \left(
          \nrcauset[3CrownBrokenLow]{5,4,2,7,1,6,3}{2/6}
          \oslash \npcauset[Singleton]{1}
        \right)
      \sqcup \left(
          \npcauset[2CrownLinked3Wedge]{7,6,5,2,8,1,4,3}
          \oslash \npcauset[Singleton]{1}
        \right)
      \sqcup \left(
          \npcauset[2Chain3Chain]{2,1,5,4,3}
          \oslash \npcauset[Singleton]{1}
        \right)
  \nexteq
  \label{eq:Example.ParallelComposition.Singleton}
    &= \nrcauset[3Crown]{4,2,6,1,5,3}{2/5}
      \sqcup \npcauset[M]{4,2,5,1,3}
      \sqcup \npcauset[2Chain]{1,2}
    \eqend{.}
  \end{align}
  \end{subequations}
\end{example}
We will come back to this example after the introduction of more general symmetries.

\subsection{Cyclic symmetries}
\label{subsec:LocalSymmetries.Cycles}

The automorphisms that act only on a finite number of elements constitute finite cycles in the automorphism group leading us to the following generalisation. 

\begin{definition}
\label{def:Symmetry.Cycles}
  Let $P$ be a poset, $Q$ be a finite poset, and $r \in \Naturals$, $r \geq 2$. 
  For an automorphism $\sigma \in \Aut(P)$, let $\Sigma(\sigma) \subseteq P$ denote the subset of all elements that are not fixed by $\sigma$. 
  A \defof{$(Q, r)$-generator} is an automorphism $\sigma \in \Aut(P)$ under the following conditions: there exists a sequence of $r$ subsets $S_i \subset \Sigma(\sigma)$, $i \in \{0, \dots, r - 1\}$ with $S_i \isom Q$, and they are the smallest, maximally ordered subsets of $\Sigma(\sigma)$ with
  \begin{align}
  \label{eq:Symmetry.SubsetCycle}
      \sigma(S_i)
    &= S_{i + 1 \mod r}
    \eqend{,}
    &
      \ell(S_i)
    &= \ell\bigl( \Sigma(\sigma) \bigr)
    \eqend{,}
    &
      \bigcup_{i = 0}^{r - 1} S_i
    &= \Sigma(\sigma)
  \end{align}
  for all $0 \leq i < r$. 
  Any distinct pair of these subsets $(S_i, S_j)$ is \defof{$(Q, r)$-symmetric}. 
  \\
  Two elements $a, b \in P$ are \defof{$(Q, r, 0)$-symmetric} if $a = b$.
  They are \defof{$(Q, r, 1)$-symmetric} if there exist subsets $A, B \subset P$ that are $(Q, r)$-symmetric with $(Q, r)$-generator $\sigma$ such that $a \in A$ and $b = \sigma^q(a) \in B$ for some $1 \leq q < r$.
  Recursively setting $n = 2, 3, \dots$, the elements are \defof{$(Q, r, n)$-symmetric} if they are not $(Q, r, j)$-symmetric for any $j < n$, but there exists some $c \in P$ and $j < n$ such that $a$ is $(Q, r, j)$-symmetric to $c$ and $c$ is $(Q, r, n - j)$-symmetric to $b$.
  \\
  For short, $a$ is \defof{$(Q, r)$-symmetric} to $b$ if there exists an $n \in \Naturals_0$ such that they are $(Q, r, n)$-symmetric. 
  For a linked subset $A \isom Q$, the \defof{$(Q, r)$-symmetry group of $A$} is the largest subgroup of $(Q, r)$-generators such that its action fixes every element $b \in P$ if there exists no element $a \in A$ that is $(Q, r)$-symmetric to $b$. 
  The \defof{$(Q, r)$-symmetry set of $A$} is the complement of all elements fixed under the action of the $(Q, r)$-symmetry group of $A \isom Q$. 
\end{definition}

Note that the condition ``smallest, maximally ordered'' in the definition of a $(Q, r)$-generator is merely a matter of convenience.
Since a representation of an automorphism in terms of a poset $Q$ and an integer $r$ is not unique in general, this condition reduces the number of possible pairs $(Q, r)$ to consider.

Because a $(Q, r)$-generator is an automorphism $\sigma \in \Aut(P)$ that fixes all elements in the complement of $\Sigma(\sigma)$, any two $(Q, r)$-symmetric elements $a, b \in \Sigma(\sigma)$ also share the same (linked) preceding and succeeding elements in $P \setminus \Sigma(\sigma)$,
\begin{subequations}
\label{eq:Symmetry.LinkSubsets}
\begin{align}
    L^\pm(a) \setminus \Sigma(\sigma)
  &= L^\pm(b) \setminus \Sigma(\sigma)
  \eqend{,}
\nexteq
    J^\pm(a) \setminus \Sigma(\sigma)
  &= J^\pm(b) \setminus \Sigma(\sigma)
  \eqend{.}
\end{align}
\end{subequations}
This is a generalisation of the singleton-symmetry, see \eqref{eq:Symmetry.Element}, or conversely, the singleton-symmetry is the special case with $Q = \npcauset[Singleton]{1}$ and $r = 2$.

The $(Q, r)$-symmetry groups are normal subgroups of automorphisms as $(Q, r)$-symmetric elements form equivalence classes. 
For a poset $(P, \leq_P)$, the retract $P \oslash_r Q$ is isomorphic to the poset quotient poset $(\tilde{P}, \leq_{\tilde{P}})$ that is given by the equivalence classes under all $(Q, r)$-symmetries such that for all $A, B \in \tilde{P}$: 
\begin{align}
\label{eq:Symmetry.QuotientOrder}
    A \leq_{\tilde{P}} B
  &\liff \exists a \in A, b \in B: 
    a \leq_P b
  \eqend{.}
\end{align}
If $r = 2$, a $(Q, 2)$-symmetry is a reflection of subsets isomorphic to $Q$, so we say $Q$-symmetric for short, and drop the index on the division operator, $P \oslash Q$. 

\begin{example}
\label{eg:ParallelComposition.Continued}
  We retract the parallel composition in \autoref{eg:ParallelComposition.Singleton} with further 2-chain- and singleton-symmetries, 
  \begin{subequations}
  \eqseqlabel{eq:Example.ParallelComposition.FurtherSymmetries}
  \begin{align}
        \left(
          \nrcauset[3Crown]{4,2,6,1,5,3}{2/5}
        \sqcup \npcauset[M]{4,2,5,1,3}
        \sqcup \npcauset[2Chain]{1,2}
        \right)
      \oslash \npcauset[2Chain]{1,2}
    &=  \npcauset[2Chain]{1,2}
        \sqcup \npcauset[Wedge]{2,1,3}
        \sqcup \npcauset[2Chain]{1,2}
    \eqend{,}
  \nexteq
        \left(
          \npcauset[2Chain]{1,2}
        \sqcup \npcauset[Wedge]{2,1,3}
        \sqcup \npcauset[2Chain]{1,2}
        \right)
      \oslash \npcauset[Singleton]{1}
    &= \left( \npcauset[2Chain_2Chain_2Chain]{5,6,3,4,1,2} \right)
    \eqend{,}
  \nexteq
        \left( \npcauset[2Chain_2Chain_2Chain]{5,6,3,4,1,2} \right)
      \oslash \npcauset[2Chain]{1,2}
    &= \npcauset[2Chain]{1,2}
    \eqend{.}
  \end{align}
  Note that the elements of $\nrcauset[3Crown]{4,2,6,1,5,3}{2/5}$ are $\npcauset[2Chain]{1,2}$-symmetric and $(\npcauset[2Chain]{1,2}, 3)$-symmetric, both with the 2-chain as retract. 
\end{subequations}
\end{example}
\begin{example}
\label{eg:3Chain2Cycle}
  The smallest poset with a 3-cycle but without a 2-cycle is 
  \begin{align}
  \label{eq:Example.3Chain2Cycle}
      \nrcauset[3Chain3Cycle]{6,8,3,1,5,9,7,2,4}{1/5,5/9,1/9/{}}
      \oslash_3 \npcauset[3Chain]{1,2,3}
    &= \npcauset[3Chain]{1,2,3}
    \eqend{,}
  \end{align}
  where the diagonal edge across the central element in the Hasse diagram links the corner elements to each other but does \emph{not} link them to the centre. 
\end{example}

From the definition, it might be obvious to a groupologist that 
\begin{proposition}
	The automorphism group of a finite poset $P$ is generated by $(Q, r)$-symmetries for a given set of $(Q, r)$ pairs.
\end{proposition}
\proof
Let $\sigma \in \Aut(P)$ be any non-trivial automorphism. 
Since any automorphism is part of a cycle (that can be represented in a cycle graph of the group), there is some integer $r \geq 2$ such that $\sigma^r = 1_{\Aut(P)}$ (at least as long as $P$ is finite). 
As in the definition, we denote the set of those elements that are not fixed by $\sigma$ as $\Sigma(\sigma) \subseteq P$. 
Consider a sequence of smallest, maximally ordered subsets $S_i \subset \Sigma(\sigma)$ such that \eqref{eq:Symmetry.SubsetCycle} holds for all $i \in \{0, 1, \dots, r - 1\}$. 
Such a sequence can be found by selecting a longest chain $S_0$ in $\Sigma(\sigma)$ and adding more elements to $S_0$, keeping it maximally ordered, until the sequence $(S_0, S_1, \dots, S_{r - 1})$ with $S_i = \sigma^i(S_0)$ covers $\Sigma(\sigma)$. 
Repeating this process for all choices of longest chains and maximally ordered extensions leads to a choice that minimises $\lvert S_0 \rvert$. 
By this construction, all of the subsets $S_i$ are isomorphic to each other and they are isomorphic to a (finite) poset $Q$. 
So the automorphism $\sigma$ is part of a $(Q, r)$-symmetry of $P$. 
For a finite poset $P$, we can do this construction for all automorphisms, so the automorphisms corresponding to $(Q, r)$-symmetries form the entire automorphism group. 
\qed

\begin{remark}
  Local symmetries can be generalised to other types of orders.
  For example, one may replace the condition of reflexivity of the ordering relation on a set $P$ by irreflexivity making $P$ a strict poset, or by more general conditions like almost reflexivity (defined in~\cite{GuoPaychaZhang:2023}).
  Because the reflexivity condition is actually not required in the definition, it can be dropped altogether making the order relation $\leq$ merely transitive and antisymmetric.
\end{remark}

\subsection{Local symmetry properties}
\label{subsec:LocalSymmetries.Properties}

To get towards a classification of posets by these symmetries, we define the following.
\begin{definition}
  A poset $P$ is \defof{$(Q, r)$-retractable} to a poset $\tilde{P}$ if $\tilde{P} = P \oslash_r Q \neq P$.
  For short, $P$ is \defof{locally $(Q, r)$-symmetric} if there exists some retract $\tilde{P}$ such that $P$ is $(Q, r)$-retractable to $\tilde{P}$, $P$ is \defof{locally symmetric} if there exists some retract $\tilde{P}$ and some pair $(Q, r)$ such that $P$ is $(Q, r)$-retractable to $\tilde{P}$, and it is \defof{locally unsymmetric} if it is not locally symmetric.
  In general, $P$ is \defof{retractable} to a poset $\tilde{P}$ if there exist some sequence of $(Q_i, r_i)$-symmetries such that $\hat{P} = P \oslash_{r_1} Q_1 \oslash_{r_2} Q_2 \oslash_{r_3} \dots \neq P$.
\end{definition}

\begin{example}
  Finite, locally unsymmetric posets have a trivial automorphism group.
  In particular, all finite chain posets (total orders) are locally unsymmetric, because there is only one element per layer.
\end{example}

\begin{remark}[Algebraic structure]
  Let $\mathfrak{P}_{[m, n]} := \left\{ P \in \mathfrak{P} \bypred m \leq \lvert P \rvert \leq n \right\}$. 
  For any $n \geq 1$, $r \geq 2$, $ \mathfrak{P}_{[1, n]}, \oslash_r)$ is a magma and the locally unsymmetric posets in $\mathfrak{P}_{[m, n]}$ with $m = \lfloor \frac{n}{r} \rfloor + 1$ are right identities. 
\end{remark}

As we want to classify (especially finite) posets by their local symmetries, it is of avail to reduce the number of relevant local symmetries, because pairs of elements can be part of different local symmetries simultaneously.
\begin{example}
\label{eg:SeriesComposition.3Chain2Chain2Chain.Quotients}
  Any series compositions of $n$ antichains (each with more than one element) is $k$-chain-retractable to the $n$-chain for any $k \leq n$, like 
  \begin{align}
  \label{eq:Example.SeriesComposition.3Chain2Chain2Chain.Quotients}
      \npcauset[3Chain2Chain2Chain]{3,2,1,5,4,7,6} \oslash \npcauset[3Chain]{1,2,3}
    = \npcauset[3Chain2Chain2Chain]{3,2,1,5,4,7,6} \oslash \npcauset[2Chain]{1,2}
    = \npcauset[3Chain2Chain2Chain]{3,2,1,5,4,7,6} \oslash \npcauset[Singleton]{1}
    &= \npcauset[3Chain]{1,2,3}
    \eqend{.}
  \end{align}
\end{example}
So we make a distinction:
\begin{definition}
  Let $P$ be a locally $(Q, r)$-symmetric poset (for some finite poset $Q$ and $r \geq 2$) and let $S \subseteq P$ be a corresponding $(Q, r)$-symmetry set (for some subset $A \subset P$). 
  The $(Q, r)$-symmetry of $S$ is \defof{composite} if there exists some finite poset $\hat{Q}$, $\hat{r} \geq 2$, and a corresponding $(\hat{Q}, \hat{r})$-symmetry set $\hat{S} \subseteq S$ such that 
  \begin{subequations}
  \label{eq:Symmetry.Composite}
  \begin{align}
      \hat{Q}
    &\subset Q
      \eqend{,}\text{~or} 
  \nexteq
      \hat{Q}
    &= Q \text{~and~} \hat{r} < r
      \eqend{.}
  \end{align}
  \end{subequations}
  It is \defof{prime} if it is not composite. 
\end{definition}

Our main focus shall lie on prime $(Q, r)$-symmetries, which implies that all $(Q, r)$-symmetric subsets are linked. 
To retract only prime $(Q, r)$-symmetries, we write a prime as in $P \oslash'_r Q$. 
In \autoref{eg:SeriesComposition.3Chain2Chain2Chain.Quotients}, the only prime retraction is by the singleton-symmetry.

As a consequence of $(Q, r)$-generators being elements of the automorphism group of a poset, the length of the longest chain subset in the poset and its $(Q, r)$-retracted version is always the same. 

\begin{proposition}[Common layers]
\label{prop:CommonLayers}
  For any poset $P$, finite poset $Q$ and $r \in \Naturals$, the symmetry quotient $P / (Q, r)$ is layer-preserving, 
  \begin{align}
  \label{eq:CommonLayers}
        \forall E \in P / (Q, r):
        \forall e \in E:
    \qquad
        \ell_{P / (Q, r)}(E)
    &= \ell_{P}(e)
    \eqend{.}
  \end{align}
\end{proposition}

This implies that any locally $(Q, r)$-symmetric poset $P$ and the retract $P \oslash_r Q$ have the same number of layers, $\ell(P) = \ell(P \oslash_r Q)$. 

For any finite, locally unsymmetric poset $Q$, the elements of $\bigsqcup_{i = 1}^{n} Q$ (with $n \geq 2$) are all pairwise $Q$-symmetric. 
More generally, the definition of $(Q, r)$-symmetric elements simplifies whenever there is some $k$ such that two $(Q, r)$-symmetric elements are at most $(Q, r, k)$-symmetric. 
As an example, consider 2-chain-symmetric elements.

\begin{lemma}
\label{lma:2ChainSymmetry}
  If two elements in a poset $P$ are $(\npcauset[2Chain]{1,2}, 2, k)$-symmetric, then $k \in \{0, 1\}$.
\end{lemma}
\proof
Let $a \equiv_1 b$ mean that $a$ is $(\npcauset[2Chain]{1,2}, 2, 1)$-symmetric to $b$. 
The only non-trivial part is to show that, for distinct elements $a, b, c \in P$, $a \equiv_1 b$ and $b \equiv_1 c$ implies $a \equiv_1 c$.

By \autoref{prop:CommonLayers}, the elements $a$, $b$ and $c$ have to be on the same layer within $P$ as well as on the same layer within the subsets isomorphic to $Q$. 
This means that they are either all on the first or all on the second layer in the 2-chain subsets.
Since the 2-chain is self-opposite, let us assume that all elements are on the first level of the 2-chain subsets, otherwise the arguments have to be made with the opposite of $P$ and all signs in the link subsets are reversed.

Let $A = \{a, \hat{a}\}$ and $B_a = \{b, \hat{b}_a\}$ such that $a \lessd \hat{a}$, $b \lessd \hat{b}_a$, $\hat{a} \neq \hat{b}_a$ and the corresponding automorphism $\alpha$ such that $\alpha(A) = B_a$.
Similarly, let $B_c = \{b, \hat{b}_c\}$ and $C = \{c, \hat{c}\}$ such that $b \lessd \hat{b}_c$, $c \lessd \hat{c}$ and let $\beta$ be the corresponding 2-chain-reflection such that $\beta(B_c) = C$.
If and only if $a < \hat{b}_a$, then $b < \hat{a}$ and the 2-chain-reflection is not prime but decomposes into two singleton-symmetries.
Any singleton-symmetry in a poset corresponds to a commuting subgroup of automorphisms that are reflections, so $k$ is at most 1 in that case.

So now we consider the case of prime symmetries $a \equiv_1 b$ and $b \equiv_1 c$.
By the definition of the 2-chain-symmetry, we have
\begin{subequations}
\label{eq:2ChainSymmetry}
\begin{align}
\label{eq:2ChainSymmetry.FirstLayer.ab}
    L^-(a)
  &= L^-(b)
  \eqend{,}
  &
    L^+(a) \setminus \{\hat{a}\}
  &= L^+(b) \setminus \{\hat{b}_a\}
  \eqend{,}
\nexteq
\label{eq:2ChainSymmetry.SecondLayer.ab}
    L^-(\hat{a}) \setminus \{a\}
  &= L^-(\hat{b}_a) \setminus \{b\}
  \eqend{,}
  &
    L^+(\hat{a})
  &= L^+(\hat{b}_a)
  \eqend{,}
\nexteq
\label{eq:2ChainSymmetry.FirstLayer.bc}
    L^-(b)
  &= L^-(c)
  \eqend{,}
  &
    L^+(b) \setminus \{\hat{b}_c\}
  &= L^+(c) \setminus \{\hat{c}\}
  \eqend{,}
\nexteq
\label{eq:2ChainSymmetry.SecondLayer.bc}
    L^-(\hat{b}_c) \setminus \{b\}
  &= L^-(\hat{c}) \setminus \{c\}
  \eqend{,}
  &
    L^+(\hat{b}_c)
  &= L^+(\hat{c})
  \eqend{.}
\end{align}
\end{subequations}
The first layer, minus equations \eqref{eq:2ChainSymmetry.FirstLayer.ab} and \eqref{eq:2ChainSymmetry.FirstLayer.bc} imply $L^-(a) = L^-(c)$. 
For the other subsets, we have to consider two cases, with the simplest examples \npcauset[2Chain_2Chain_2Chain]{5,6,3,4,1,2} (separated case) and \nrcauset[3Crown]{4,2,6,1,5,3}{2/5} (crossing case), respectively.

In the separated case, $\hat{b}_a = \hat{b}_c$, the subsets in \eqref{eq:2ChainSymmetry.FirstLayer.ab} and \eqref{eq:2ChainSymmetry.FirstLayer.bc}, as well as \eqref{eq:2ChainSymmetry.SecondLayer.ab} and \eqref{eq:2ChainSymmetry.SecondLayer.bc} are pairwise identical. 
So $a \equiv_1 c$ given the 2-chain-symmetric subsets $A = \{a, \hat{a}\}$ with $a \lessd \hat{a}$ and $C = \{c, \hat{c}\}$ with $c \lessd \hat{c}$. 

In the crossing case, $\hat{b}_a \neq \hat{b}_c$, \eqref{eq:2ChainSymmetry.FirstLayer.ab}$^+$ implies $\hat{b}_c \in L^+(a) \setminus \{\hat{a}\}$ so $\hat{b}_c \notin \{\hat{a}, \hat{c}\}$, and \eqref{eq:2ChainSymmetry.FirstLayer.bc}$^+$ implies $\hat{b}_a \in L^+(c) \setminus \{\hat{c}\}$ so $\hat{b}_a \notin \{\hat{a}, \hat{c}\}$.
Then \eqref{eq:2ChainSymmetry.SecondLayer.ab}$^-$ implies $c \in L^-(\hat{a})$ so also $\hat{a} \in L^+(c)$, as well as \eqref{eq:2ChainSymmetry.SecondLayer.bc}$^-$ implies $a \in L^-(\hat{c})$ so also $\hat{c} \in L^+(a)$.
We find that $\hat{a} = \hat{c}$, because $\hat{a}, \hat{c} \notin L^+(b)$ by the first layer identities \eqref{eq:2ChainSymmetry.FirstLayer.ab}$^+$ and \eqref{eq:2ChainSymmetry.FirstLayer.bc}$^+$, and $\hat{a} \neq \hat{c}$ contradicts these equations in combination.
Combining \eqref{eq:2ChainSymmetry.FirstLayer.ab}$^+$ and \eqref{eq:2ChainSymmetry.FirstLayer.bc}$^+$ together with the union by $\{\hat{a}\} = \{\hat{c}\}$, we get
\begin{subequations}
\begin{align}
    \Bigl( L^+(a) \setminus \{\hat{a}, \hat{b}_c\} \Bigr) \cup \{\hat{a}\}
  &= \Bigl( L^+(c) \setminus \{\hat{b}_a, \hat{c}\} \Bigr) \cup \{\hat{c}\}
\nexteq
    L^+(a) \setminus \{\hat{b}_c\}
  &= L^+(c) \setminus \{\hat{b}_a\}
  \eqend{.}
\end{align}
\end{subequations}
Similarly, \eqref{eq:2ChainSymmetry.SecondLayer.ab}$^-$ and \eqref{eq:2ChainSymmetry.SecondLayer.bc}$^-$ combine to give
\begin{subequations}
\begin{align}
    \Bigl( L^-(\hat{b}_a) \setminus \{b, c\} \Bigr) \cup \{b\}
  &= \Bigl( L^-(\hat{b}_c) \setminus \{a, b\} \Bigr) \cup \{b\}
\nexteq
    L^-(\hat{b}_a) \setminus \{c\}
  &= L^-(\hat{b}_c) \setminus \{a\}
  \eqend{,}
\nexteq
    L^+(\hat{b}_a)
  &= L^+(\hat{b}_c)
  \eqend{.}
\end{align}
\end{subequations}
With $A = \{a, \hat{b}_c\}$, $C = \{c, \hat{b}_a\}$ as the 2-chains $a \lessd \hat{b}_c$ and $c \lessd \hat{b}_a$ that are mapped to each other by a 2-chain-reflection, we obtain $a \equiv_1 c$.
\qed

\subsection{Classes of posets with local symmetries}
\label{subsec:Classification}

All posets that are $(Q, r)$-retractable to some poset $R$ form a class of \emph{symmetry extensions} (or \emph{symmetry class})
\begin{align}
\label{eq:PosetClass.1}
    [R \odot_r Q]
  &:= \left\{
      P \in \mathfrak{P}
    \bypred
      P \oslash_r Q = R \neq P
    \right\}
  \eqend{,}
\end{align}
and we write $[R \odot_r Q]'$ to include only prime symmetries.
These posets are extensions of $R$ by (prime) $(Q, r)$-symmetries.
In the following, we consider only prime symmetries.

Larger classes are obtained with sequences of symmetry extensions, like
\begin{subequations}
\label{eq:PosetClass.N}
\begin{align}
\label{eq:PosetClass.2}
    [R \odot_{r_1} Q_1 \odot_{r_2} Q_2]'
  := \bigl\{
      P \in \mathfrak{P}
    \bigm|
      P \oslash'_{r_2} Q_2
    &= R_1 \neq P \eqend{,} \nnexteq
      R_1 \oslash'_{r_1} Q_1
    &= R \neq R_1
    \bigr\}
  \eqend{,}
\nexteq
\label{eq:PosetClass.3}
    [R \odot_{r_1} Q_1 \odot_{r_2} Q_2 \odot_{r_3} Q_3]'
  := \bigl\{
      P \in \mathfrak{P}
    \bigm|
      P \oslash'_{r_3} Q_3
    &= R_2 \neq P \eqend{,} \nnexteq
      R_2 \oslash'_{r_2} Q_2
    &= R_1 \neq R_2 \eqend{,} \nnexteq
      R_1 \oslash'_{r_1} Q_1
    &= R \neq R_1
    \bigr\}
  \eqend{,}
\end{align}
\end{subequations}
and so on to arbitrary order.

\begin{example}
\label{eg:AntichainClass}
  The simplest symmetry class collects all antichains that are singleton-retractable to the singleton, see \autoref{eg:Antichains.Singleton}, 
  \begin{align}
  \label{eq:PosetClass.1Chain}
      [\npcauset[Singleton]{1} \odot \npcauset[Singleton]{1}]'
    &= \Bigl\{
        (\npcauset[2Antichain]{2,1}),
        (\npcauset[3Antichain]{3,2,1}),
        (\npcauset[4Antichain]{4,3,2,1}),
        (\npcauset[5Antichain]{5,4,3,2,1}),
        \dots
      \Bigr\}
    \eqend{.}
  \end{align}
  Since this symmetry class already includes all 1-layer posets (except for the locally unsymmetric singleton) and the layers are preserved by retractions there are no other, non-trivial classes of 1-layer posets.
  In particular, 
  \begin{align}
  \label{eq:PosetClass.1Chain1Chain}
      [\npcauset[Singleton]{1}
        \odot \npcauset[Singleton]{1}
        \odot \npcauset[Singleton]{1}
      ]'
    &= \emptyset
    \eqend{.}
  \end{align}
\end{example}

More generally, for any poset $R$ that is not singleton-retractable, the class $[R \odot \npcauset[Singleton]{1}]'$ includes all extensions by the singleton-symmetry and 
\begin{align}
\label{eq:PosetClass.RepeatedElementSymmetry}
    [R \odot \npcauset[Singleton]{1} \odot \npcauset[Singleton]{1}]'
  &= \emptyset
  \eqend{.}
\end{align}
The singleton-symmetry is always prime since it is the smallest symmetry (and the prime symbol is not really necessary, but kept to emphasise the focus on prime symmetries here). 
Similarly to the antichains symmetry class, 
\begin{lemma}[Repeated parallel composition]
\label{lma:RepeatedParallelComposition}
  For every locally unsymmetric poset $R$, 
  \begin{align}
  \label{eq:RepeatedParallelComposition}
      \left\{
        \bigsqcup_{i = 1}^{n} R
      \bypred
        n \geq 2
      \right\}
    \subseteq [R \odot R]'
    \eqend{.}
  \end{align}
\end{lemma}

The order in the sequence of symmetry retractions in \eqref{eq:PosetClass.N} is important, because symmetry retractions do not commute and there are posets that are retractable to more than one locally unsymmetric poset. 
For example,
\begin{subequations}
\begin{align}
    \nrcauset[3CrownBrokenLow]{5,4,2,7,1,6,3}{2/6}
    \oslash \npcauset[Singleton]{1}
    \oslash \npcauset[2Chain]{1,2}
  &= \nrcauset[3Crown]{4,2,6,1,5,3}{2/5}
    \oslash \npcauset[2Chain]{1,2}
  = \npcauset[2Chain]{1,2}
  \eqend{,}
\nexteq
    \nrcauset[3CrownBrokenLow]{5,4,2,7,1,6,3}{2/6}
    \oslash \npcauset[2Chain]{1,2}
    \oslash \npcauset[Singleton]{1}
  &= \npcauset[2ChainLinkedWedge]{4,3,1,5,2}
    \oslash \npcauset[Singleton]{1}
  = \npcauset[N]{2,4,1,3}
  \eqend{.}
\end{align}
\end{subequations}
This is because, the singleton-symmetry in the poset on the left partially ``breaks'' the 2-chain-symmetry of $\nrcauset[3Crown]{4,2,6,1,5,3}{2/5}$. 
The sequences of symmetry retractions can also be different but still yield the same poset, as in
\begin{subequations}
\begin{align}
    \nrcauset[Double2ChainLinkedDoubleVee]{9,4,8,2,7,10,5,3,1,6}{2/5,4/7}
    \oslash \npcauset[3Wedge]{3,2,1,4}
    \oslash \npcauset[2Chain]{1,2}
  &= \npcauset[DownhookBrokenDouble2Chain]{5,3,4,1,2,6}
    \oslash \npcauset[2Chain]{1,2}
  = \npcauset[Downhook]{3,1,2,4}
  \eqend{,}
\nexteq
    \nrcauset[Double2ChainLinkedDoubleVee]{9,4,8,2,7,10,5,3,1,6}{2/5,4/7}
    \oslash \npcauset[Vee]{1,3,2}
    \oslash \npcauset[Wedge]{2,1,3}
  &= \npcauset[Double2ChainLinkedVee]{6,2,5,7,3,1,4}
    \oslash \npcauset[Wedge]{2,1,3}
  = \npcauset[Downhook]{3,1,2,4}
  \eqend{.}
\end{align}
\end{subequations}
More so, the retraction sequences do not even have to have the same length, see the example~\eqref{eq:KROrders.ExampleContractions} in a later section.
This demonstrates that the set of all posets $\mathfrak{P}$ does not partition into locally unsymmetric posets and classes formed from symmetry extensions of these.
Nevertheless, individual poset classes are often easier to enumerate, and as long as we keep track of repeated countings, we can use poset classes to enumerate all posets for a fixed cardinality $n$.

Define the counting function as
\begin{align}
\label{eq:CountingFunction}
    c_n : \mathcal{P}(\mathfrak{P})
  &\to \Naturals_0
  \eqend{,}
&
    \mathfrak{C}
  &\mapsto \Bigl|
      \bigl\{
        P \in \mathfrak{C}
      \bigm\vert
        | P | = n
      \bigr\}
    \Bigr|
  \eqend{.}
\end{align}
For \autoref{eg:AntichainClass}, we have
\begin{align}
\label{eq:CountingFunction.1Chain.1Chain}
    c_n[\npcauset[Singleton]{1} \odot \npcauset[Singleton]{1}]'
  &= \begin{cases}
      0 & \text{if~} n \in \{0, 1\} \eqend{,} \\
      1 & \text{otherwise} \eqend{.}
    \end{cases}
\end{align}
For any finite, locally unsymmetric poset $R$ with cardinality $m := \lvert R \rvert$, the counting function of the symmetry class under the singleton-symmetry vanishes for $n \leq m$ and, for $n > m$,
\begin{align}
\label{eq:CountingFunction.LocallyUnsymmetricPoset.1Chain}
    c_n[R \odot \npcauset[Singleton]{1}]'
  &= \binom{(n - m) + m - 1}{n - m}
  = \binom{n - 1}{m - 1}
  \eqend{,}
\end{align}
because there are $n - m$ elements to be distributed into the $m$ equivalence classes formed by the elements of the poset $R$ (while the order in which the elements are distributed does not matter).

\section{Posets prime retractable to the 2-chain}
\label{sec:2ChainClass}

In this section, we explore posets (prime) retractable to the 2-chain that form the simplest symmetry classes other than $[\npcauset[Singleton]{1} \odot \npcauset[Singleton]{1}]'$.

\subsection{Complete bipartite graphs and crown posets}
\label{subsec:2ChainClass.CompleteBipartitesAndCrowns}

The simplest symmetry class of posets prime retractable to the 2-chain are the posets that correspond to complete bipartite graphs, for which all elements of the first and second layer form an equivalence class by the singleton-symmetry, 
\begin{align}
\label{eq:PosetClass.2Chain.1Chain}
    [\npcauset[2Chain]{1,2} \odot \npcauset[Singleton]{1}]'
  &= \left\{
      \npcauset[Wedge]{2,1,3},
      \npcauset[Vee]{1,3,2},
      \npcauset[3Wedge]{3,2,1,4},
      \npcauset[2Gonal]{2,1,4,3},
      \npcauset[3Vee]{1,4,3,2},
      \npcauset[4Wedge]{4,3,2,1,5},
      \npcauset[3Chain2Chain]{3,2,1,5,4},
      \npcauset[2Chain3Chain]{2,1,5,4,3},
      \npcauset[4Vee]{1,5,4,3,2},
      \dots
    \right\}
  \eqend{.}
\end{align}
The enumeration of this symmetry class is a special case of \eqref{eq:CountingFunction.LocallyUnsymmetricPoset.1Chain}, so for $n > 2$, 
\begin{align}
\label{eq:CountingFunction.2Chain.1Chain}
    c_n[\npcauset[2Chain]{1,2} \odot \npcauset[Singleton]{1}]'
  &= n - 1
  \eqend{.}
\end{align}
The enumeration of all bipartite graphs was considered in detail in~\cite{Hanlon:1979}. 

By \autoref{lma:RepeatedParallelComposition}, the parallel compositions of 2-chains are 2-chain-retractable, 
\begin{align}
\label{eq:PosetClass.2Chain.2Chain.ParallelCompositionSubclass}
    [\npcauset[2Chain]{1,2} \odot \npcauset[2Chain]{1,2}]'
  &\supseteq \Bigl\{
      \left( \npcauset[2Chain_2Chain]{3,4,1,2} \right),
      \left( \npcauset[2Chain_2Chain_2Chain]{5,6,3,4,1,2} \right),
      \left( \npcauset[2Chain_2Chain_2Chain_2Chain]{7,8,5,6,3,4,1,2} \right),
      \dots
    \Bigr\}
  \eqend{.}
\end{align}
Furthermore, the prime symmetry class $[\npcauset[2Chain]{1,2} \odot \npcauset[2Chain]{1,2}]'$ is complemented by another family of posets. 

\begin{definition}[Crown posets]
\label{def:CrownPosets}
  The \defof{$n$-crown poset} $X_n$ ($n \geq 2$) is the poset with two layers of elements $a_1, a_2, \dots, a_n$ and $b_1, b_2, \ldots, b_n$ and the partial order $\leq$ such that 
  \begin{align}
  \label{eq:CrownPosets}
      \forall i, j
    &\in [1, n]:
  &
      i \neq j
    &\liff a_i < b_j
    \eqend{.}
  \end{align}
\end{definition}

Any crown poset corresponds to a complete bipartite graph where the edges between opposite vertices have been removed. 
Note that the $n$-crown poset (for $n > 2$) contains the $k$-crown posets for all $k < n$. 

\begin{lemma}[Crown retractions]
\label{lma:CrownPosetRetractions}
  The $n$-crown poset $X_n$ is 2-chain-retractable to the 2-chain. 
\end{lemma}
\proof
By definition of the $n$-crown, the event $a_i$ is not linked to $b_i$ for all $i \in [1, n]$. 
Considering the two 2-chains $(a_i, b_j)$ and $(a_j, b_i)$ for any $i \neq j$, the element pairs $(a_i, a_j)$ and $(b_i, b_j)$ are 2-chain symmetric, because $a_i$ and $a_j$ are both linked to all $b_k$ for $k \notin \{i, j\}$, and similarly $b_i$ and $b_j$ are both linked to all $a_k$ for $k \notin \{i, j\}$. 
This holds for all $i \neq j$ so that all elements $a_i$ form a single equivalence class, and so do all $b_i$.
\qed

\begin{theorem}
\label{thm:PosetClass.2Chain.2Chain}
  \begin{subequations}
  \label{eq:PosetClass.2Chain.2Chain}
  \begin{align}
      [\npcauset[2Chain]{1,2} \odot \npcauset[2Chain]{1,2}]'
    &= \left\{
        \bigsqcup_{i = 1}^{n} \npcauset[2Chain]{1,2}
      \bypred
        n \geq 2
      \right\}
      \cup \left\{ X_n \bypred n \geq 3 \right\}
  \nexteq
    &= \left\{
        \left( \npcauset[2Chain_2Chain]{3,4,1,2} \right),
        \left( \npcauset[2Chain_2Chain_2Chain]{5,6,3,4,1,2} \right),
        \left( \npcauset[2Chain_2Chain_2Chain_2Chain]{7,8,5,6,3,4,1,2} \right),
        \dots,
        \nrcauset[3Crown]{4,2,6,1,5,3}{2/5}, 
        \nrcauset[4Crown]{5,3,2,8,1,7,6,4}{2/7,3/6},
        \nrcauset[5Crown]{6,4,3,2,10,1,9,8,7,5}{2/9,3/8,4/7},
        \nrcauset[6Crown]{7,5,4,3,2,12,1,11,10,9,8,6}{2/11,3/10,4/9,5/8},
        \dots
      \right\}
    \eqend{.}
  \end{align}
  \end{subequations}
\end{theorem}
\proof
For any poset $P$ to be prime 2-chain-retractable to the 2-chain, it has to have two layers by \autoref{prop:CommonLayers} and all elements in the first and second layer each have to form a single equivalence class in $P / (\npcauset[2Chain]{1,2}, 2)$ with all elements in the equivalence classes being $(\npcauset[2Chain]{1,2}, 2, 0)$-symmetric or $(\npcauset[2Chain]{1,2}, 2, 1)$-symmetric by \autoref{lma:2ChainSymmetry}. 
This is true for the repeated parallel composition, see \autoref{lma:RepeatedParallelComposition}, and it is true for every crown by \autoref{lma:CrownPosetRetractions}. 
It remains to show that there is no other poset in the symmetry class. 
We will show that any poset $P$ in this symmetry class that is not a parallel composition of 2-chains is a crown poset. 

Since $P$ is not a parallel composition of 2-chains, it has to be connected and it must contain two 2-chain subsets such that the second/first layer elements are linked by an element in the first/second layer, so $\npcauset[M]{4,2,5,1,3}$ or $\npcauset[W]{3,5,1,4,2}$ is a subset of $P$, respectively. 
Under the quotient by the 2-chain, these subsets form two equivalence classes on the first or second layer, respectively, and one equivalence class on the other layer. 
Only if $P$ has yet another element on the second resp.\ first layer that is linked to both elements on the other layer, the quotient results in one equivalence class per layer. 
In either case, this means the 3-crown is a subset of $P$. 

The remaining argument is recursive assuming that we have an $n$-crown as a subset of $P$ for some $n \geq 3$. 
Note that if $P$ is in the symmetry class, so is its opposite, because the 2-chain is self-opposite. 
We denote the elements on the first layer $a_1, a_2, \dots, a_n$ and the elements on the second layer $b_1, b_2, \dots, b_n$ (as in \autoref{def:CrownPosets}). 
If $P$ contains other elements, there has to be another one, say $a_{n + 1}$ linked to the $n$-crown because $P$ is connected, and we assume it is on the first layer (without loss of generality). 
To keep the first order 2-chain symmetry of all elements, $a_{n + 1}$ has to be linked to every element $b_i$ with $i \in [1, n]$. 
As before, the resulting subset now has two equivalence classes on the first layer when taking the quotient by the 2-chain. 
Only if there is also another element on the second layer $b_{n + 1}$ linked to every element $a_i$ with $i \in [1, n]$ but not linked to $a_{n + 1}$, the quotient is isomorphic to the 2-chain. 
The subset formed by $a_1, a_2, \dots, a_{n + 1}$, $b_1, b_2, \dots, b_{n + 1}$ and their partial order is the $(n + 1)$-crown poset. 
We can repeat this argument until all elements of $P$ are included, so $P$ is a crown poset. 
Since all crown posets are already included, the symmetry class is complete. 
\qed

For any cardinality $n$, there are 
\begin{align}
\label{eq:CountingFunction.2Chain.2Chain}
    c_n[\npcauset[2Chain]{1,2} \odot \npcauset[2Chain]{1,2}]'
  &= \begin{cases}
      2 & \text{if $n$ is even and~} n \geq 6 \eqend{,} \\
      1 & \text{if~} n = 4 \eqend{,} \\
      0 & \text{otherwise}
    \end{cases}
\end{align}
posets that are prime 2-chain retractable to the 2-chain. 

To get larger symmetry classes (prime) retractable to the 2-chain with further singleton-symmetries and 2-chain-reflections, recall that any repeated extension by singleton-symmetries does not yield any posets, see~\eqref{eq:PosetClass.RepeatedElementSymmetry}, but we may have alternating symmetry quotients as in \autoref{eg:ParallelComposition.Continued}, or repeat 2-chain-extensions.
\begin{example}
\label{eg:Repeated2ChainSymmetry}
  The following poset retracts to the 2-chain by a sequence of repeated 2-chain quotients, 
  \begin{align}
      \nrcauset[3CrownBroken3Crown]{7,5,3,10,2,9,1,8,6,4}{1/6,2/4,3/8,5/9}
      \oslash \npcauset[2Chain]{1,2}
      \oslash \npcauset[2Chain]{1,2}
    &= \nrcauset[3Crown]{4,2,6,1,5,3}{2/5}
        \oslash \npcauset[2Chain]{1,2}
    = \npcauset[2Chain]{1,2}
    \eqend{.}
  \end{align}
\end{example}

Next, we consider posets that are $(Q, r)$-symmetric with $Q = \npcauset[2Chain]{1,2}$ and $r \geq 3$.

\subsection{Fences and polygons}
\label{subsec:2ChainClass.FencesAndPolygons}

As any poset that is retractable to the 2-chain has to have two layers (by \autoref{prop:CommonLayers}), we analyse another type of connected, 2-layer posets.

\begin{theorem}[Symmetries of fences]
\label{thm:Fences}
  Let 
  \begin{align}
      f : \Naturals
    &\to \Naturals
    \eqend{,}
    &
      n
    &\mapsto \begin{cases}
        n & \text{if $n$ is even}, \\
        \frac{n + 1}{2} & \text{if $n$ is odd},
      \end{cases}
  \end{align}
  and let $f^{\downarrow}$ be the repeated composition $f \after f \after f \after \cdots$ until the result is even.
  Every fence of even cardinality is locally unsymmetric.
  Every fence of odd cardinality $n$ is retractable to the $f^{\downarrow}(n)$-fence.
\end{theorem}
\proof
The $n$-fence $F_n$ with even $n \geq 4$ has exactly two elements that have only one link each, $a_1$ and $a_n$, that are on different layers so that any automorphism $\alpha \in \Aut(F_n)$ has to leave these fixed.
So we could as well remove these elements from the fence to obtain the $(n - 2)$-fence without changing the automorphism group, $\Aut(F_{n - 2}) \isom \Aut(F_n)$.
As this holds for any even $n$ and the smallest fence, the 2-chain, is locally unsymmetric, the automorphism group of all even fences is trivial.

For odd $n = 2 \hat{n} + 1$, the $n$-fences $F_n^{\vee}$ and $F_n^{\wedge}$ have two elements $a_1$ and $a_n$ on the same layer.
The fences $F_3^{\vee}$ and $F_3^{\wedge}$ are singleton-symmetric and retract to the 2-chain, which is $F_2$, and $f(n) = 2$.
For $\hat{n} \geq 2$, let $A$ be the $\hat{n}$-fence subset of $F_n^{\vee}$ ($F_n^{\wedge}$) including $a_1$ and $B$ be the $\hat{n}$-fence subset of $F_n^{\vee}$ ($F_n^{\wedge}$) including $a_n$.
The retraction by $F_{\hat{n}}^{\vee}$ ($F_{\hat{n}}^{\wedge}$) folds and matches the subset $B$ onto $A$, while keeping the central element $a_{\hat{n} + 1}$ fixed, and $F_n \oslash Q$ is the upper (lower) $f(n)$-fence.
If $f(n)$ is odd, repeat the retraction to get an $f^2(n)$-fence, and so on, until the resulting fence has even cardinality $f^{\downarrow}(n)$ and is locally unsymmetric by the first statement. 
\qed

\autoref{thm:Fences} implies that, for any $n \in \Naturals$, the upper and the lower $(2^n + 1)$-fences are retractable to the 2-chain.
Except for the 3-fences (\npcauset[Wedge]{2,1,3} and \npcauset[Vee]{1,3,2}), the retractions require a sequence of more than one quotient.
Though there are other 2-layer posets that have fence-symmetries and retract to the 2-chain by only \emph{one} quotient.
In particular:
\begin{definition}
\label{def:PolygonalPosets}
  The \defof{$n$-gonal poset} $G_n$ is the $2 n$-fence with an additional link that makes it cyclic,
  \begin{align}
  \label{eq:PolygonalPoset}
    a_1 < a_2 > a_3 < \dots a_{2 n}
  &> a_1
    \eqend{.}
  \end{align}
\end{definition}

Note that in the literature you may find these polygonal posets referred to as crowns, in contrast to the generalisation of the 3-gonal poset that I call crown posets in \autoref{def:CrownPosets}.
The definition of polygonal posets is to highlight the abstract representations of regular geometrical objects in terms of posets below.
The $n$-gonal poset represents the hull of a regular $n$-gon and a series composition with a singleton on the third layer yields the closure of the polygon that includes the central 2-face.
It is obvious from the definition of the $n$-gonal poset that it is retractable to the 2-chain by a $(\npcauset[2Chain]{1,2}, n)$-symmetry.
Additionally, the automorphism groups of a polygon also admits subgroups that are generated by reflections only.

\begin{theorem}[Reflection symmetries of polygonal posets]
\label{thm:PolygonalPosets}
  The $n$-gonal poset is retractable to the 2-chain by an $(n - 1)$-fence-reflection (followed by a retraction of the singleton-symmetry if $n$ is even).
\end{theorem}
\begin{figure}
  \centering
  \includegraphics{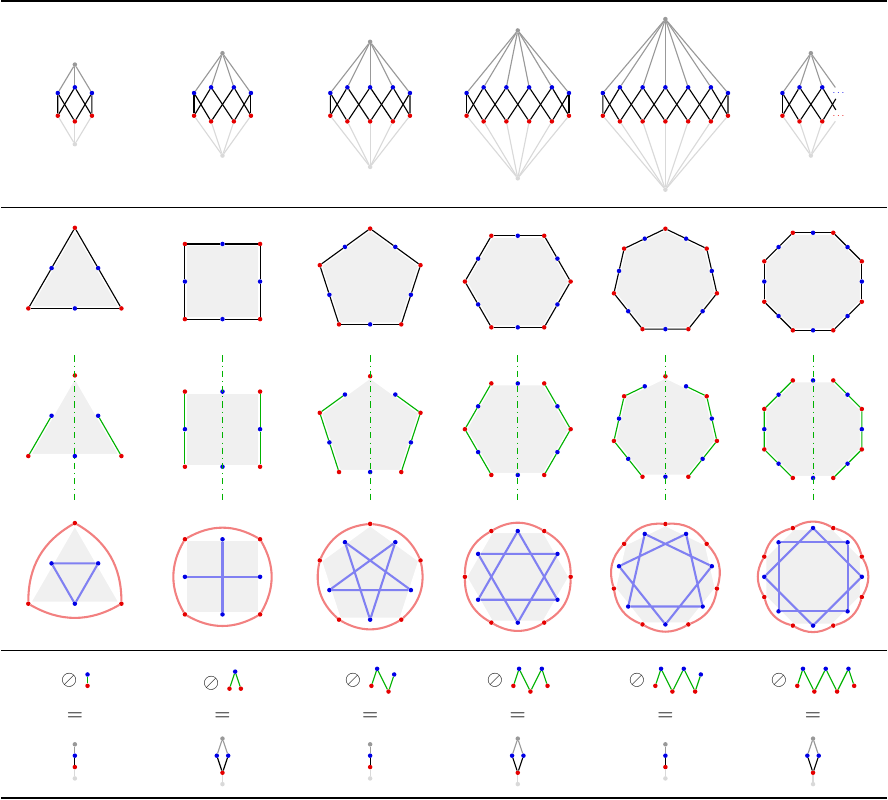}
  \caption{\label{fig:PolygonPosets} The first-layer and second-layer elements of polygonal posets (first row) correspond to the vertices (red) and edges (blue) of regular polygons (second row), respectively. 
  Regular polygons have reflection symmetries (dash dotted axis) that are represented by lower fences (green, third row). 
  Taking the quotient by the fence is equivalent to matching vertex and edge elements for all reflection axes that meet an edge. 
  This results in one equivalence class (fourth row) for the vertices (red cycle) and one or two equivalence classes --- the conjugacy classes of the respective dihedral groups --- represented by the non-degenerate or degenerate star polygons (blue), respectively. 
  Hence, the quotients are isomorphic to the 2-chain or the Vee-poset (fifth row).
  The additional element in each poset (first row) above the polygonal poset is the closure including the central 2-face (grey).
  The element below each poset (light grey) is the dual to the top face, it is the abstract $(-1)$- or least face (empty shape).}
\end{figure}
\proof
The $n$-gonal poset $G_n$ corresponds to the vertices (first layer elements) and edges (second layer elements) of a regular polygon, see \autoref{fig:PolygonPosets}. 
Let $A \subset G_n$ such that $A \isom F_{n - 1}^{\wedge}$. 
There are two elements $c_1$ and $c_2$ that are linked to some elements in $A$, which fall on the symmetry axis (dash dotted line) in the third row of \autoref{fig:PolygonPosets}.
Note that $B := G_n \setminus A \setminus \{c_1, c_2\}$ is also isomorphic to $F_{n - 1}^{\wedge}$.
The subsets $A$ and $B$ are symmetric under reflection about the symmetry axis as shown in \autoref{fig:PolygonPosets}.
There are reflected subsets $A$ and $B$ for any symmetry axis that meets an edge element.
All such reflections form a subgroup of the dihedral group $\mathrm{D}_n$ (symmetry group of the $n$-gon).
The subgroup is the entire dihedral group $\mathrm{D}_n$ if $n$ is odd and it is half the dihedral group if $n$ is even.
The retraction corresponds to reducing the $n$-gon symmetries to the center of the dihedral group $\operatorname{Z}(\mathrm{D}_n)$, which is only one element (the identity) if $n$ is odd, $G_n \oslash F_{n - 1}^{\wedge} = \npcauset[2Chain]{1,2}$, and two elements if $n$ is even, $G_n \oslash F_{n - 1}^{\wedge} = \npcauset[Vee]{1,3,2} \in [\npcauset[2Chain]{1,2} \odot \npcauset[Singleton]{1}]'$.
\qed

Taking the quotient by the subgroup of the dihedral group is equivalent to removing the inner automorphisms because $\mathrm{D}_n / \operatorname{Z}(\mathrm{D}_n) \isom \operatorname{Inn}(\mathrm{D}_n)$ \cite{Miller:1942}.
Under the action of the dihedral (sub)group, any vertex of the $n$-gon is mapped to any other vertex, while any edge is mapped to any other edge only in the odd case. 
This is graphically represented by connecting the poset elements of the $n$-gon edges to a star $n$-gon that splits into two $\frac{n}{2}$-gons if $n$ is even, see fourth row in \autoref{fig:PolygonPosets}.
The resulting equivalence classes in $G_n / (Q, 2)$ are the conjugacy classes of the dihedral group.
Similarly, one may take the quotients with upper fences, where the axes of symmetry have to meet at least one vertex element resulting in a 2-chain or the poset $\npcauset[Wedge]{2,1,3} \in [\npcauset[2Chain]{1,2} \odot \npcauset[Singleton]{1}]'$, respectively.

Each polygonal poset is self-opposite, see first row in \autoref{fig:PolygonPosets}. 
Also when including the element for the 2-face as the closure of the geometric shape (grey, above) as well as its dual, least element (light grey, below), the resulting poset is self-opposite.
This is equivalent to the fact that the geometric shape of a regular polygon is self-dual, and each vertex can be uniquely matched with an edge.

The polygon theorem generalises to those polyhedra that are derived from the infinite family of polygons and have dihedral symmetry, too: the self-dual pyramids, bipyramids dual to prisms, and antiprisms dual to trapezohedrons.
However, posets corresponding to the vertex-edge pairing are not retractable to the 2-chain (but to the N-poset) in general.
Exceptions are those polygons that have larger symmetry groups, like the tetrahedron, the 3-cube or the octahedron.

\subsection{Regular polytopes and dimensional reduction}
\label{subsec:2ChainClass.Polytopes}

\begin{figure}
  \centering
  \includegraphics{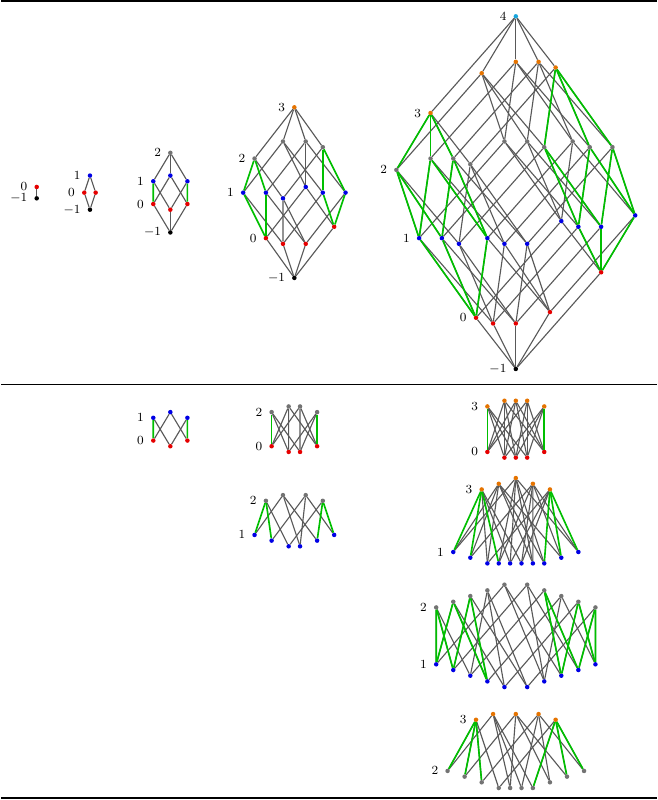}
  \caption{\label{fig:SimplexPosets} The simplices in dimensions 0 to 4 as posets (first row), where the dimension index of the faces is placed to the left of each layer.
  The rows below show all (but opposite) distinct subsets formed by the elements of two distinct layers.
  Green links mark examples of subsets that are reflected.}
\end{figure}

In any dimension $d \geq 3$, there is a finite number of (convex) regular polytopes.
In $d = 3$, we have the five Platonic solids and in $d = 4$, there are six polytopes including the three special cases of the 24-cell, 120-cell and the 600-cell.
There are three families of polytopes that exists in all dimensions $d \geq 3$, the self-dual simplices, and the hypercubes dual to the orthoplexes.
\begin{definition}[Simplex posets]
  For any dimension $d \in \Naturals_0$, the \defof{$d$-simplex poset} $S_d$ is the poset of $d + 2$ layers with 1 element on the first layer, the least face, preceding $d + 1$ elements on the second layer for the vertices of the $d$-simplex. 
  The third layer contains an element for each edge of the $d$-simplex such that it succeeds exactly those elements from the first layer that correspond to the boundary (the two vertices) of the edge. 
  This continues to each layer $k$ with an element for each $(k - 2)$-face succeeding the elements corresponding to its boundary. 
\end{definition}
The Hasse diagram of the $d$-simplex poset is self-opposite and partitions into two copies of the $(d - 1)$-simplex poset that are offset by 1 layer, giving the lower left and upper right subsets for each of the smallest simplices in the first row of \autoref{fig:SimplexPosets}.

\begin{theorem}[Reflection symmetries of simplex posets]
\label{thm:SimplexPosets}
  For any $d \geq 2$, the $d$-simplex poset is $(d - 2)$-simplex-retractable to the $(d + 2)$-chain.
\end{theorem}
\proof
For $d \geq 2$, the $d$-simplex has a reflection symmetry about any hyperplane that contains any of the $(d - 2)$-faces (together with all its boundary faces, its closure).
The reflection matches the two adjacent $(d - 2)$-faces and all their boundary faces (except for those meeting the symmetry plane), which are represented by $(d - 2)$-simplex subsets (green in the first row of \autoref{fig:SimplexPosets}).
Because the choice of the hyperplane is otherwise arbitrary, there is a matching for any pair of faces of the same dimension.
Thus the corresponding quotient by the $(d - 2)$-simplex poset yields a retraction to the $(d + 2)$-chain (for each face type, there is one equivalence class).
\qed

\begin{corollary}
  For the $d$-simplex poset, each subset of two layers is retractable to the 2-chain. 
\end{corollary}
\proof
For subsets of all elements from two distinct layers of the $d$-simplex, taking the quotient by the respective subsets of the quotient poset yields posets that retract to the 2-chain. 
\qed

\begin{example}
  The 2-layer subsets of $d$-simplexes with $d \in \{2, 3, 4\}$ are shown in the lower rows of \autoref{fig:PolygonPosets}.
  For each simplex, there is a pair of reflected subposets marked in green.
  We may take any subposet of all elements in a layer (corresponding to all faces of the same type).
  For example, the 2-layer subset of all vertices and all $(d - 1)$-faces (second row in \autoref{fig:SimplexPosets}) has, for each vertex, only one non-neighbouring $(d - 1)$-face, so that the subset is the $(d + 1)$-crown poset.
\end{example}

The retraction of the $d$-simplex to a chain can be understood as a reduction from a higher dimensional order to a $1$-dimensional (or total) order, where the order dimension is the number of total orders that intersect to the partial order, see \cite{DushnikMiller:1941}.
For the application to causal set theory further below, we may consider the dimension of a partial order to be given in terms of an embedding in Minkowski spacetime.
\begin{definition}
\label{def:MinkowskiEmbeddingDimension}
Let $\Minkowski^d$ be $d$-dimensional Minkowski space ($\Reals^{1, d - 1}$ with Lorentzian metric tensor $\eta = \diag(1, -1, -1, -1, \dots)$, which is equipped with a partial ordering determined by its causal structure.
The \defof{Minkowski embedding dimension} of a poset $P$ is the smallest value $d \in \Naturals$ such that there exists an order-preserving embedding $P \to \Minkowski^d$. 
\end{definition}
For example, the singleton poset embeds in 0-dimensional Minkowski space, any larger chain embeds in 1-dimensional Minkowski space (a timeline), and a $d$-simplex poset is embeddable in $(1 + d)$ dimensions.
Note that the poset representation has a dimension one larger than the dimension of the geometric objects, as the extra dimension distinguishes different types of faces.

\begin{figure}
  \centering
  \includegraphics{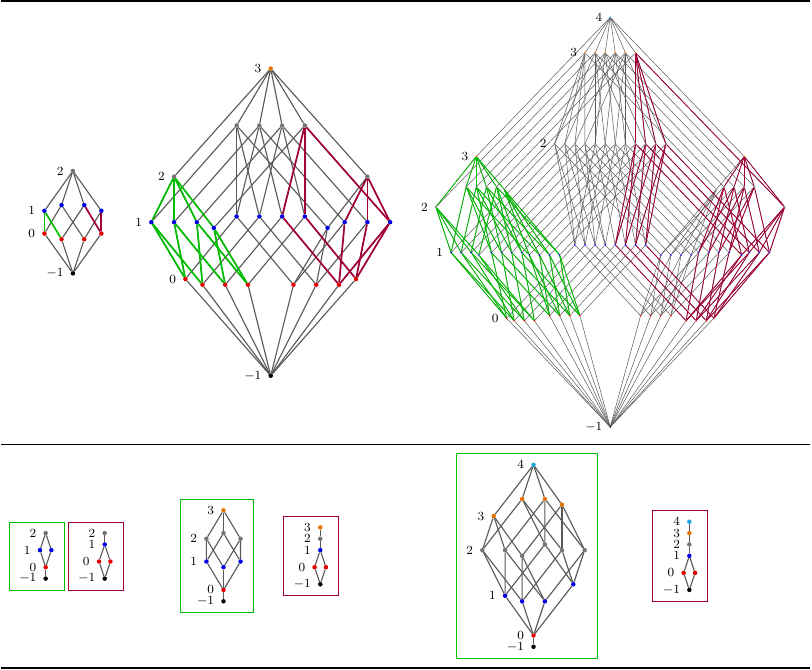}
  \caption{\label{fig:CubePosets} The square, 3-cube and 4-cube represented by their Hasse diagrams (first row), where the dimension of the faces is indicated to the left of each layer.
  The second row shows the posets that are obtained when taking a quotient by two different types of reflection (hyper)planes.
  Different to \autoref{fig:SimplexPosets}, here only one of subposet is shown for each of the two reflections.
  The green subposet is reflected by to a hyperplane that cuts through edges, while the purple subposet is reflected by a hyperplane that meets vertices.}
\end{figure}
Cubes and orthoplexes are the other two families of regular polytopes that exist in any dimension $d \geq 3$.
Since orthoplexes are the dual polytopes to cubes (they have opposite posets), results about their local symmetries are also dual to the local symmetries of the cubes.

The $d$-cube poset is defined completely analogous to the $d$-simplex.
Three examples are shown in \autoref{fig:CubePosets}, where two different subsets for symmetry retractions correspond to the green and purple subposets, respectively.
Taking a quotient by either of these subsets yields the two retracted posets in the second row.

For each $d$-cube, the retract shown on the left (green) is the series composition of the singleton with the $(d - 1)$-simplex, $\npcauset[Singleton]{1} \vee S_{d - 1}$.
Geometrically, the retraction by this subposet correspond to the symmetry reflections about hyperplanes that pass through edges (and higher faces) but not through vertices.
Any vertex of the $d$-cube is mapped to any other vertex by a sequence of such reflections.
While it takes only one reflection to map any vertex to any of its adjacent vertices, opposite pairs of vertices require $d$ reflections.

The retract shown on the right (magenta) is the series composition of the singleton with the 2-antichain and a $d$-chain, see \autoref{fig:CubePosets}.
Here any subset of two layers that does not include the vertex layer is retractable to the 2-chain by a single symmetry quotient.
In terms of local symmetries (\autoref{def:Symmetry.Cycles}), the $d$-cube is a principal example for elements that are $(Q, 2, d)$-symmetric, where $Q$ is isomorphic to the magenta selection on the right and the geometric dimension $d$ is the order of the symmetry relation.

\subsection{Local symmetries by complete bipartite graphs}
\label{subsec:2ChainClass.CompleteBipartite}

The 2-layer poset composed of the 1- and 2-face-elements of the 3-simplex poset $S_3$, see \autoref{fig:SimplexPosets}, is prime wedge-retractable to the 2-chain.
\begin{theorem}
\label{thm:PosetClass.2Chain.Wedge}
  \begin{align}
  \label{eq:PosetClass.2Chain.Wedge}
      \left[ \npcauset[2Chain]{1,2} \odot \npcauset[Wedge]{2,1,3} \right]'
    &= \left\{ \nrcauset[3Simplex12Faces]{8,6,4,3,10,2,9,1,7,5}{2/7,3/7,3/9,4/5,4/10,6/9} \right\}
    \eqend{.}
  \end{align}
\end{theorem}
\proof
Similarly to the proof of \autoref{thm:PosetClass.2Chain.2Chain}, we only need to show that there are no other posets in the symmetry class.
Any poset $P$ in the symmetry class has to be connected, because the only non-connected posets that are wedge-retractable to the 2-chain are parallel compositions of wedges, but these are also singleton-retractable, so their wedge-symmetries are not prime.
Since the connected poset $P$ is wedge-retractable to the 2-chain, it has to contain a subset of two wedges that are linked by another element, for which there are three options:
\begin{align}
\label{eq:PosetClass.2Chain.Wedge.ConnectedOptions}
    O_1
  &:= \npcauset[7Fence]{6,4,7,2,5,1,3}
  \eqend{,}
&
    O_2
  &:= \npcauset[3WedgeLinkedWedge]{6,5,3,7,2,1,4}
  \eqend{,}
&
    O_3
  &:= \npcauset[2GonalLinkedWedge]{5,4,7,2,1,6,3}
  \eqend{.}
\end{align}
More so, all second layer elements have to be first order wedge-symmetric to retract to the 2-chain by the wedge-symmetry, because the wedge contains a 2-chain as a subposet, has only one element on the second layer, so that \autoref{lma:2ChainSymmetry} applies.
For the first option $O_1$ (lower 7-fence), the two elements from the left on the second layer are first order wedge-symmetric if there are two additional elements on the first layer such that the resulting subset is
\begin{align}
\label{eq:PosetClass.2Chain.Wedge.ConnectedOption.Extended}
    O_1'
  &:= \nrcauset[3Simplex12FacesIncomplete]{8,6,4,3,9,2,7,1,5}{3/5,3/9,4/7}
  \eqend{.}
\end{align}

In the second option $O_2$, both second layer elements are already wedge-symmetric.
But to have the second and third element from left on the first layer wedge-symmetric, we also require an additional element on the second layer, leading again to $O_1'$. 

The wedge-symmetry between the first and second element from the left on the second layer in the third option $O_3$ requires two additional elements on the first layer that are both linked to the first and third element on the second layer (from the left). 
This extended subset is the first poset in \eqref{eq:PosetClass.2Chain.Wedge.ConnectedOption3.Tiaras} below, and any further enlargement follows a similar argument as in \autoref{thm:PosetClass.2Chain.2Chain} leading to a sequence of crowns where each element on the first layer is replaced by a 2-antichain: 
\begin{align}
\label{eq:PosetClass.2Chain.Wedge.ConnectedOption3.Tiaras}
    \nrcauset[321Tiara]{7,6,4,3,9,2,1,8,5}{3/8,4/8},
    \nrcauset[421Tiara]{9,8,6,5,4,3,12,2,1,11,10,7}{3/11,4/11,5/10,6/10},
    \nrcauset[521Tiara]{11,10,8,7,...,3,15,2,1,14,13,12,9}{3/14,4/14,5/13,6/13,7/12,8/12},
    \dots
  \eqend{.}
\end{align}
However, these subsets do not have a prime wedge-symmetry, which renders the third option irrelevant. 

The subset $O_1'$ completes to the edges--2-faces poset of the 3-simplex 
\begin{align}
\label{eq:PosetClass.2Chain.Wedge.Edge2Faces}
    S_3|_{1,2}
  &:= \nrcauset[3Simplex12Faces]{8,6,4,3,10,2,9,1,7,5}{2/7,3/7,3/9,4/5,4/10,6/9}
\end{align}
by requiring an element on the second layer linked to the three elements on the first layer with a single link in $O_1'$.
This is the only choice to make all elements on the first layer wedge-symmetric and there cannot be any other elements in the subset, so that $P = S_3|_{1,2}$.

If there where another element $x$ on the second layer, it would have to be linked to all first-layer elements of $S_3|_{1,2}$ to maintain the wedge-symmetry of these.
So $x$ and each other second-layer element $y$ share three of their succeeded elements on the first layer of $S_3|_{1,2}$.
Let $a_1, a_2, a_3$ be the other three elements linked to $x$ but not to $y$.
To have $x$ first-order wedge-symmetric to $y$, we can at most take two of these elements to form a wedge subset with $x$ at the tip, say $a_1$ and $a_2$.
The third element $a_3$ would also have to be linked to $y$ to make $y$ wedge-symmetric to $x$, but this is not the case in the subset $S_3|_{1,2}$.
So we reached a contradiction.

Similarly, an additional element on the first layer would have to be linked to all second-layer elements of $S_3|_{1,2}$, so that there must be another second-layer element as well.
So $S_3|_{1,2}$ is the only connected poset $P$ that is prime wedge-retractable to the 2-chain.
\qed

\begin{corollary}
\label{cor:PosetClass.2Chain.Vee}
  As the vee poset \npcauset[Vee]{1,3,2} is opposite to the wedge poset, the theorem implies
  \begin{align}
  \label{eq:PosetClass.2Chain.Vee}
      \left[ \npcauset[2Chain]{1,2} \odot \npcauset[Vee]{1,3,2} \right]'
    &= \left\{ \nrcauset[3Simplex01Faces]{6,4,10,2,9,1,8,7,5,3}{1/7,2/5,2/8,4/8,4/9,6/7} \right\}
    \eqend{.}
  \end{align}
\end{corollary}

Obviously, their counting functions are 
\begin{align}
\label{eq:CountingFunction.2Chain.WedgeVee}
		c_n\left[ \npcauset[2Chain]{1,2} \odot \npcauset[Wedge]{2,1,3} \right]'
	= c_n\left[ \npcauset[2Chain]{1,2} \odot \npcauset[Vee]{1,3,2} \right]'
	&= \begin{cases}
			1 & \text{if~} n = 10 \eqend{,} \\
			0 & \text{otherwise} \eqend{.}
		\end{cases}
\end{align}

Similar symmetry classes result by series composing an antichain $A_n$ with $n \geq 2$ with a singleton, $Q = A_n \vee \npcauset[Singleton]{1}$.
The $Q$-symmetry (or $\op{Q}$-symmetry) extensions of the 2-chain are given by corresponding 2-layer subposets shown in \autoref{fig:SimplexPosets}. 
For example, the symmetry class of prime 3-wedge-retractable posets contains two posets of the 4-simplex subsets formed from the 1/3-faces and 2/3-faces, 
\begin{align}
\label{eq:PosetClass.2Chain.3Wedge}
    \left[ \npcauset[2Chain]{1,2} \odot \npcauset[3Wedge]{3,2,1,4} \right]'
  &= \left\{
      \ncauset[4Simplex23Faces]{13,11,9,7,6,5,4,15,3,14,2,12,1,10,8}{1/8,1/10,2/8,2/12,3/8,3/14,4/8,4/15,5/10,5/12,6/10,6/14,7/10,7/15,9/12,9/14,11/12,11/15,13/14,13/15},
      \nrcauset[4Simplex13Faces]{12,10,8,7,6,5,4,3,15,2,14,1,13,11,9}{2/13,3/13,3/14,4/11,4/15,5/11,5/14,6/11,6/13,7/9,7/15,8/9,8/14,10/13}
    \right\}
  \eqend{.}
\end{align}
With similar arguments as in the proof of \autoref{thm:PosetClass.2Chain.Wedge}, it can be shown that there are no \emph{prime} extensions of the 2-chain by $Q$-symmetries where $Q = A_m \vee A_n$ for some antichains $A_m$, $A_n$ (with $m, n \geq 2$). 

Using local symmetries for the classification of posets reduces the problem of enumerating all posets to finding and enumerating all locally unsymmetric posets (and then building symmetry classes from these).
Finding all locally unsymmetric posets for a given cardinality becomes essential.
In the following, we list locally unsymmetric posets and find their enumeration (at least for small cardinalities).

\section{Enumeration of locally unsymmetric posets by layers}
\label{sec:Enumeration}

In this section, we compare the numbers of locally unsymmetric posets to all posets with a fixed cardinality, grouped by the number of layers.
We first consider the distribution of numbers of posets as function of their number of layers and then determine the exact enumeration of posets with a number of layers close to their cardinality.
In the end, there are some remarks on the enumeration of posets with few layers, since these dominate for large cardinalities.

\subsection{Numbers of locally unsymmetric posets vs.\ all posets}
\label{sec:Enumeration.LocallyUnsymmetricVsAll}

Let $\mathfrak{U}$ denote the set of all locally unsymmetric posets. 
The smallest of these posets (up to cardinality $n = 5$) are 
\begin{subequations}
\label{eq:LocallyUnsymmetricPosets}
\begin{align}
    \mathfrak{U}_{1}
  &= \{\npcauset[Singleton]{1}\}
  \eqend{,}
\nexteq
    \mathfrak{U}_{2}
  &= \{\npcauset[2Chain]{1,2}\}
  \eqend{,}
\nexteq
    \mathfrak{U}_{3}
  &= \Bigl\{
      \left( \npcauset[2Chain_Singleton]{3,1,2} \right),
      \npcauset[3Chain]{1,2,3}
    \Bigr\}
  \eqend{,}
\nexteq
    \mathfrak{U}_{4}
  &= \biggl\{
      \npcauset[N]{2,4,1,3},
      \left( \npcauset[3Chain_Singleton]{4,1,2,3} \right),
      \npcauset[Downhook]{3,1,2,4},
      \npcauset[Uphook]{1,4,2,3},
      \npcauset[4Chain]{1,2,3,4}
    \biggr\}
  \eqend{,}
\nexteq
    \mathfrak{U}_{5}
  &\pickindent{=
    \Biggl\{}
      \left( \npcauset[N_Singleton]{5,2,4,1,3} \right),
      \left( \npcauset[Downhook_Singleton]{5,3,1,2,4} \right),
      \left( \npcauset[Uphook_Singleton]{5,1,4,2,3} \right),
      \left( \npcauset[3Chain_2Chain]{4,5,1,2,3} \right),
      \npcauset[NUphook]{4,1,5,2,3},
      \npcauset[NDownhook]{3,5,1,2,4},
      \npcauset[Hook]{4,1,3,5,2},
      \npcauset[Ndown]{4,1,2,5,3},
      \npcauset[Nup]{2,5,1,3,4},
      \npcauset[Fish]{3,1,5,2,4},
      \npcauset[2Diamonddescend]{3,1,4,2,5},
      \npcauset[2Diamondascend]{1,3,5,2,4},
  \eqbreakr
      \left( \npcauset[4Chain_Singleton]{5,1,2,3,4} \right),
      \npcauset[Downhookdown]{4,1,2,3,5},
      \npcauset[Uphookup]{1,5,2,3,4},
      \npcauset[Downhookup]{3,1,2,4,5},
      \npcauset[Uphookdown]{1,2,5,3,4},
      \npcauset[2Chaindiamond]{1,4,2,3,5},
      \npcauset[5Chain]{1,2,3,4,5}
    \Biggr\}
  \eqend{.}
\end{align}
\end{subequations}

For the enumeration of locally unsymmetric posets among all posets $\mathfrak{P}$, we define the following functions and numbers. 
The number of posets for a given cardinality $n$ and number of layers $l$ is 
\begin{align}
\label{eq:PosetNumbers.CardinalityLayer}
      p(n, l)
  &:= \Bigl\lvert
        \bigl\{
          P \in \mathfrak{P}
        \bigm\vert
          \lvert P \rvert = n, 
          \ell(P) = l
        \bigr\}
      \Bigr\rvert
  \eqend{.}
\end{align}
All posets with fixed cardinality but an arbitrary number of layers is the integer sequence A000112 in the OEIS~\cite{OEIS:A000112:2023}, 
\begin{align}
\label{eq:PosetNumbers.Cardinality}
      p_{n}
  &:= \Bigl\lvert
        \bigl\{
          P \in \mathfrak{P}
        \bigm\vert
          \lvert P \rvert = n
        \bigr\}
      \Bigr\rvert
  = \sum_{l = 1}^{n}
        p(n, l)
  \eqend{.}
\end{align}
Similarly, the numbers of locally unsymmetric posets $\mathfrak{U}$ are defined as 
\begin{subequations}
\label{eq:PosetNumbers.LocallyUnsymmetric}
\begin{align}
      u(n, l)
  &:= \Bigl\lvert
        \bigl\{
          P \in \mathfrak{U}
        \bigm\vert
          \lvert P \rvert = n, 
          \ell(P) = l
        \bigr\}
      \Bigr\rvert
  \eqend{,}
\nexteq
      u_{n}
  &:= \Bigl\lvert
        \bigl\{
          P \in \mathfrak{U}
        \bigm\vert
          \lvert P \rvert = n
        \bigr\}
      \Bigr\rvert
  = \sum_{l = 1}^{n}
        u(n, l)
  \eqend{.}
\end{align}
\end{subequations}

\begin{table}
  \centering
  \begin{subtable}{0.95\textwidth}
    \centering
    \begin{tabular}{r|*{7}r|r}
    \toprule
    $p(n, l)$
       &    1 &    2 &     3 &     4 &     5 &     6 &     7 & $p_n$ \\  
    \midrule
     1 &    1 &      &       &       &       &       &       &     1 \\  
     2 &    1 &    1 &       &       &       &       &       &     2 \\  %
     3 &    1 &    3 &     1 &       &       &       &       &     5 \\  %
     4 &    1 &    8 &     6 &     1 &       &       &       &    16 \\  %
     5 &    1 &   20 &    31 &    10 &     1 &       &       &    63 \\  %
     6 &    1 &   55 &   162 &    84 &    15 &     1 &       &   318 \\  %
     7 &    1 &  163 &   940 &   734 &   185 &    21 &     1 &  2045 \\  %
    \bottomrule
    \end{tabular}
    \caption{\label{tab:PosetNumbers.All} Numbers of (unlabelled) posets by layers and in total (last column).}
  \end{subtable}
  \\[0.5cm]
  \begin{subtable}{0.95\textwidth}
    \centering
    \begin{tabular}{r|*{7}r|r}
    \toprule
    $u(n, l)$
         &    1 &    2 &     3 &     4 &     5 &     6 &     7 & $u_n$ \\  
    \midrule
     1 &    1 &      &       &       &       &       &       &     1 \\  
     2 &    0 &    1 &       &       &       &       &       &     1 \\
     3 &    0 &    1 &     1 &       &       &       &       &     2 \\
     4 &    0 &    1 &     3 &     1 &       &       &       &     5 \\
     5 &    0 &    1 &    11 &     6 &     1 &       &       &    19 \\
     6 &    0 &    3 &    47 &    41 &    10 &     1 &       &   102 \\
     7 &    0 &    9 &   267 &   332 &   106 &    15 &     1 &   730 \\
    \bottomrule
    \end{tabular}
    \caption{\label{tab:PosetNumbers.LocallyUnsymmetric} Numbers of locally unsymmetric (automorphically trivial) posets by layers and in total.}
  \end{subtable}
  \caption{\label{tab:PosetNumbers} Number of posets (a) and locally unsymmetric posets (b) by cardinality $n$ (increasing along the rows) and layer count $l$ (increasing along the columns).}
\end{table}

\begin{figure}
  \centering
  \begin{subfigure}{0.95\textwidth}
    \centering
    \includegraphics{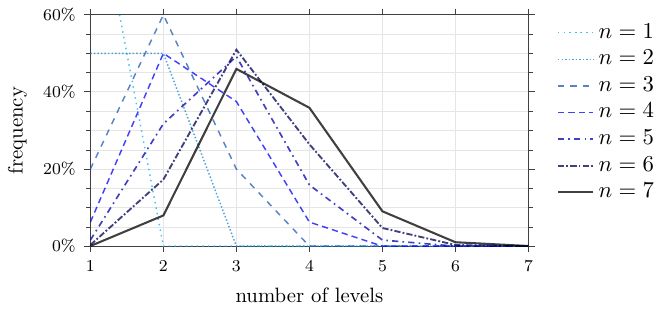}
    \caption{\label{fig:PosetNumberDistributions.All} Frequency $p n, l) / p_{n}$ vs.\ the number of layers $l$ for different cardinalities $n$.}
    \vspace{1em}
  \end{subfigure}
  \begin{subfigure}{0.95\textwidth}
    \centering
    \includegraphics{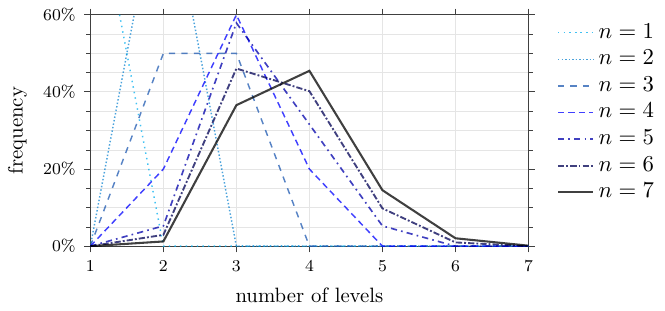}
    \caption{\label{fig:PosetNumberDistributions.LocallyUnsymmetric} Frequency $u(n, l) / u_{n}$ vs.\ the number of layers $l$ for different cardinalities $n$.}
    \vspace{0.5em}
  \end{subfigure}
  \caption{\label{fig:PosetNumberDistributions} Distributions of the number of (locally unsymmetric) posets relative to the total number of (locally unsymmetric) posets vs.\ the number of layers (length of longest chain) in the posets, all for different cardinalities $n$ of the posets.}
\end{figure}

Explicit values for all numbers up to cardinality $n = 7$ are listed in \autoref{tab:PosetNumbers} and plotted in \autoref{fig:PosetNumberDistributions}. 
Note that the peak of the distribution (as functions of the layer number $l$) for each cardinality $n \leq 7$ has its maximum at or below $l = 3$, while the peaks of the corresponding distributions for locally unsymmetric posets with the same cardinality are at higher $l$. 

\begin{remark}
  There is an equivalence between posets and coloured graphs. 
  By assigning a unique colour to each layer in any poset $P$, from color 1 to color $\ell(P)$, it is easy to see that there is a bijection between the posets (or their Hasse diagrams) and certain colourings of simple graphs. 
  The poset number $p(n, l)$ is equivalent to the number of $l$-coloured, simple graphs where any node with colour $k$ has to be graph connected to at least one node with colour $k - 1$ if $k > 1$. 
\end{remark}
\begin{example}
  There are three 2-coloured graphs with 3 nodes ($n = 3, l = 2$); the union of the single node (colour 1) and a graph-connected node pair (colour 1, colour 2); the 3-chain of 3 nodes (colour 1, colour 2, colour 1); as well as the chain of 3 nodes (colour 2, colour 1, colour 2). 
  These correspond to the posets \npcauset[2Chain_Singleton]{3,1,2}, \npcauset[Wedge]{2,1,3}, and \npcauset[Vee]{1,3,2}, respectively. 
\end{example}

For any cardinality $n$, there exists one antichain ($l = 1$),
\begin{align}
\label{eq:PosetNumbers.1Layer}
    p(n, 1)
  &= 1
  \eqend{.}
\end{align}
The singleton is locally unsymmetric and all other antichains are retractable to it,
\begin{align}
\label{eq:PosetNumbers.LocallyUnsymmetric.1Layer}
    u(n, 1)
  &= \begin{cases}
      1 & \text{if~} n = 1 \eqend{,} \\
      0 & \text{otherwise} \eqend{.}
    \end{cases}
\end{align}
The largest possible layer number is $l = n$ that corresponds to the $n$-chain, which is locally unsymmetric,
\begin{align}
\label{eq:PosetNumbers.MaxLayers}
    p(n, n)
  = u(n, n)
  &= 1
  \eqend{.}
\end{align}

\subsection{Posets with one layer less than their cardinality}
\label{subsec:Enumeration.MaxBut1Layers}

For the enumeration of posets with $l = n - 1$ layers, we find:
\begin{theorem}
\label{thm:PosetNumbers.MaxBut1Layers}
  \begin{subequations}
  \label{eq:PosetNumbers.MaxBut1Layers}
  \begin{align}
        p(n, n - 1)
    &= \frac{1}{2} n (n - 1)
    \eqend{,}
  \nexteq
        u(n, n - 1)
    &= \frac{1}{2} (n - 1) (n - 2)
    \eqend{.}
  \end{align}
  \end{subequations}
\end{theorem}
\proof
All posets with cardinality $n$ and layer number $n - 1$ are $(n - 1)$-chains with an extra element that can be separated from the chain, linked once, or linked twice.
Denote an arbitrary (possibly empty) chain by a single white vertex, so that there are the posets
\begin{align}
\label{eq:Posets.MaxBut1Layers}
    \underbrace{
      \left( \npcauset[SingletonWild_Singleton]{3,1,2/wild} \right)
    }_{n \geq 2},
    \underbrace{
      \npcauset[WedgeWild]{2,1,3,4/wild},
      \npcauset[WildVee]{1/wild,2,4,3}
    }_{n \geq 3},
    \underbrace{
      \npcauset[WildUphook]{1,5,2,3,4/wild},
      \npcauset[WildDownhook]{4,1/wild,2,3,5},
      \npcauset[WildDiamond]{1,4,2,3/wild,5}
    }_{n \geq 4},
    \underbrace{
      \npcauset[SingletonWildUphook]{1,2,6,3,4,5/wild},
      \npcauset[WildDownhookSingleton]{4,1/wild,2,3,5,6},
      \npcauset[SingletonWildDiamond]{1,2,5,3,4/wild,6},
      \npcauset[WildDiamondSingleton]{1,4,2,3/wild,5,6}
    }_{n \geq 5},
    \underbrace{
      \npcauset[2ChainWildUphook]{1,2,3,7,4,5,6/wild},
      \npcauset[WildDownhook2Chain]{4,1/wild,2,3,5,6,7},
      \npcauset[2ChainWildDiamond]{1,2,3,6,4,5/wild,7},
      \npcauset[SingletonWildDiamondSingleton]{1,2,5,3,4/wild,6,7},
      \npcauset[WildDiamond2Chain]{1,4,2,3/wild,5,6,7}
    }_{n \geq 6},
    \dots
\end{align}
where the first brace exists for $n \geq 2$ with an $(n - 2)$-chain replacing the white vertex, the second brace exists for $n \geq 3$ with an $(n - 3)$-chain replacing the white vertex, and so on.
The $k$-th brace exists for $n \geq k + 1$ and it encloses $k$ posets.
Hence, we have
\begin{align}
\label{eq:PosetNumbers.MaxBut1Layers.Derivation}
    p(n, n - 1)
  &= \sum_{k = 2}^{n} (k - 1)
  = \frac{1}{2} n (n - 1)
  \eqend{.}
\end{align}
Because there is only one element additional to the $(n - 1)$-chain subset, there are $n - 1$ posets that are retractable to the $(n - 1)$-chain whenever the additional element shares an equivalence class by a singleton-symmetry with one of the elements in the chain:
\begin{align}
\label{eq:Posets.MaxBut1Layers.LocallySymmetric}
    \underbrace{\npcauset[2AntichainWild]{2,1,3/wild}}_{n \geq 2},
    \underbrace{\npcauset[VeeWild]{1,3,2,4/wild}}_{n \geq 3},
    \underbrace{\npcauset[2Chain2AntichainWild]{1,2,4,3,5/wild}}_{n \geq 4},
    \underbrace{\npcauset[3Chain2AntichainWild]{1,2,3,5,4,6/wild}}_{n \geq 5},
    \dots
  \eqend{.}
\end{align}
Since the $(n - 1)$-chain is locally unsymmetric, we can at most find a singleton-symmetry with one more element, so there are no other local symmetries, and 
\begin{align}
\label{eq:PosetNumbers.MaxBut1Layers.LocallyUnsymmetric}
    u(n, n - 1)
  &= p(n, n - 1) - (n - 1)
  = \frac{1}{2} (n - 1) (n - 2)
  \eqend{.}
  \qed
\end{align}

The ratio of locally unsymmetric posets to all posets (for $n \geq 2$) is
\begin{align}
\label{eq:PosetNumbers.MaxBut1Layers.Ratio}
    \frac{u(n, n - 1)}{p(n, n - 1)}
  &= 1 - \frac{2}{n}
\end{align}
and approaches 1 as $n \to \infty$. 
Asymptotically, almost all $n$-element posets with $n - 1$ layers are locally unsymmetric.

\subsection{Posets with two layers less than their cardinality}
\label{subsec:Enumeration.MaxBut2Layers}

The enumeration of posets with $l = n - 2$ layers follows similarly.
\begin{theorem}
  For $n < 3$, the poset numbers $p(n, n - 2)$ and $u(n, n - 2)$ vanish, and for $n \geq 3$, they are 
  \begin{subequations}
  \label{eq:PosetNumbers.MaxBut2Layers}
  \begin{align}
      p(n, n - 2)
    &= \frac{1}{6} (n - 2)
      \Bigl( n^3 - 2 n^2 - 5 n + 12 \Bigr)
    \eqend{,}
  \nexteq
      u(n, n - 2)
    &= \frac{1}{6} (n - 3)
      \Bigl( n^3 - 4 n^2 + 2 n - 2 \Bigr)
    \eqend{.}
  \end{align}
  \end{subequations}
\end{theorem}
\proof
All $(n - 2)$-layer posets with cardinality $n$ contain an $(n - 2)$-chain and two additional elements that are -- similarly to the $(n - 1)$-layer posets -- each either separated from the chain or linked to the chain with one or two links, while the two elements can be but do not have to be linked to each other as well. 
Say that an element is linked up/down to the chain at position $i \in [1, n - 2]$ if it precedes/succeeds the event at $i$. 
For $n = 3$, there is only the 3-antichain, while for $n \geq 4$, we have the following cases 
\begin{enumerate}
\item For any poset $P$ from the sequence \eqref{eq:Posets.MaxBut1Layers} with cardinality $n - 1$, we have the parallel composition with the singleton, so there are 
  \begin{align}
  \label{eq:PosetNumbers.MaxBut2Layers.NLinks0Links}
      P \sqcup \npcauset[Singleton]{1}:
  &&
      p(n - 1, n - 2)
    &= \frac{1}{2} n^2 - \frac{3}{2} n + 1
  \end{align}
  possible combinations. 
\item One of the additional elements is linked up and the other one is linked down, or both are linked up/down,   
  \begin{align}
  \label{eq:PosetNumbers.MaxBut2Layers.1Link1Link}
      \npcauset[ChainUD]{5,1/wild,2,4/wild,6,7/wild,3},
      \underbrace{
        \npcauset[ChainDD]{6,3,1/wild,2,4/left label=$i$,5/wild,7/left label=$j$,8/wild},
        \npcauset[ChainUU]{1/wild,2,8,3/wild,4,7,5,6/wild}
      }:
  &&
      (n - 3)^2
      + 2 \sum_{i = 2}^{n - 2}
      \sum_{j = i}^{n - 2} 1
    &= 2 n^2 - 11 n + 15
    \eqend{.}
  \end{align}
\item One element is linked to the chain twice and the other element is linked to the chain once, independently, 
  \begin{align}
  \label{eq:PosetNumbers.MaxBut2Layers.2Links1Link}
      \underbrace{
        \npcauset[ChainDB]{2/wild,3/left label=$i$,7,4/wild,1,5/left label=$k$,6/wild,8/left label=$j$,9/wild},
        \npcauset[ChainUB]{1/wild,2,7,3/wild,4,6/wild,8,9/wild,5}
      }:
  &&
      2 \sum_{i = 1}^{n - 4}
      \sum_{j = i + 2}^{n - 2}
      \sum_{k = 2}^{n - 2} 1
    &= n^3 - 10 n^2 + 33 n - 36
    \eqend{.}
  \end{align}
\item Both elements are linked to the chain twice, independently, 
  \begin{subequations}
  \eqseqlabel{eq:PosetNumbers.MaxBut2Layers.2Links2Links}
  \begin{align}
      \npcauset[ChainBBCommonlow]{1/wild,2/left label=$i$,6,4,5/wild,7/left label=$j$,8/wild,3,9/left label=$k$,10/wild}:
  &&
      \sum_{i = 1}^{n - 4}
      \sum_{j = i + 2}^{n - 2}
      \sum_{k = j}^{n - 2} 1
    &= \frac{1}{6} n^3 - \frac{3}{2} n^2 + \frac{13}{3} n - 4
    \eqend{,} \text{~and~}
  \nexteq
      \npcauset[ChainBB]{1/wild,2/left label=$i$,7,3/wild,4/left label=$k$,6/wild,8/left label=$j$,9/wild,5,10/left label=$l$,11/wild}:
  &&
      \sum_{i = 1}^{n - 4}
      \sum_{j = i + 2}^{n - 2}
      \sum_{k = i + 1}^{n - 4}
      \sum_{l = k + 2}^{n - 2} 1
    &= \frac{1}{8} n^4 - \frac{23}{12} n^3 + \frac{87}{8} n^2 - \frac{325}{12} n + 25
    \eqend{.}
  \end{align}
  \end{subequations}
\item The elements are linked together and at most one of them is linked to the chain once (ignoring the posets already counted with point 2), 
  \begin{align}
  \label{eq:PosetNumbers.MaxBut2Layers.Linked01Links0Links}
      \left( \npcauset[Chain2Chain]{4,5,1,2,3/wild} \right),
      \underbrace{
        \npcauset[Chain2ChainupperD]{6,1/wild,2/left label=$i$,7,3,4,5/wild},
        \npcauset[Chain2ChainlowerU]{4,7,1/wild,2,3,5,6/wild}
      },
      \underbrace{
        \npcauset[Chain2ChainupperU]{4,5,1/wild,2,3,6/left label=$i$,7/wild},
        \npcauset[Chain2ChainlowerD]{1/wild,2,6,7,3,4,5/wild}
      }:
  &&
      1
      + 2 \sum_{i = 1}^{n - 4} 1
      + 2 \sum_{i = 3}^{n - 2} 1
    &= 4 n - 15
    \eqend{.}
  \end{align}
\item The elements are linked together and one of them is linked to the chain twice  (ignoring the posets already counted with point 3), 
  \begin{align}
  \label{eq:PosetNumbers.MaxBut2Layers.Linked0Links2Links}
      \underbrace{
        \npcauset[Chain2ChainupperB]{6,1/wild,2/left label=$i$,7,3,4,5/wild,8/left label=$j$,9/wild},
        \npcauset[Chain2ChainlowerD]{1/wild,2,6,9,3,4,5/wild,7,8/wild}
      }:
  &&
      2 \sum_{i = 1}^{n - 5}
      \sum_{j = i + 3}^{n - 2} 1
    &= n^2 - 9 n + 20
    \eqend{.}
  \end{align}
\item The elements are linked together and each of them is linked to the chain once (ignoring the posets already counted with point 3), 
  \begin{subequations}
  \eqseqlabel{eq:PosetNumbers.MaxBut2Layers.Linked1Link1Link}
  \begin{align}
      \underbrace{
        \npcauset[Chain2ChainlowerDupperD]{1/wild,2/left label=$i$,8,3/wild,4/left label=$j$,9,5,6,7/wild},
        \npcauset[Chain2ChainlowerUupperU]{4,7,1/wild,2,3,5,6/wild,8,9/wild}
      }:
  &&
      2 \sum_{i = 1}^{n - 5} \sum_{j = i + 1}^{n - 4} 1
    &= n^2 - 9 n + 20
    \eqend{,} \text{~and~}
  \nexteq
      \npcauset[Chain2ChainlowerDupperU]{1/wild,2/left label=$i$,6,7,3,4/wild,5,8/left label=$j$,9/wild},
      \npcauset[Chain2ChainlowerUupperD]{4,1/wild,2/left label=$i$,7,3/wild,5/left label=$j$,6/wild}:
  &&
      \sum_{i = 1}^{n - 5}
      \sum_{j = i + 3}^{n - 2} 1
      + \sum_{i = 1}^{n - 3}
      \sum_{j = i + 1}^{n - 2} 1
    &= n^2 - 7 n + 13
    \eqend{.}
  \end{align}
  \end{subequations}
\item The elements are linked together, one is linked to the chain twice, and the other is linked to the chain once (ignoring the posets already counted with point 4), 
  \begin{subequations}
  \eqseqlabel{eq:PosetNumbers.MaxBut2Layers.Linked1Link2Links}
  \begin{align}
      \underbrace{
        \npcauset[Chain2ChainlowerBupperD]{1/wild,2/left label=$i$,6,3/wild,4/left label=$j$,9,5/wild,7/left label=$k$,8/wild},
        \npcauset[Chain2ChainlowerUupperB]{4,1/wild,2,7,3/wild,5,6/wild,8,9/wild}
      }:
    &&
      2 \sum_{i = 1}^{n - 4}
      \sum_{j = i + 1}^{n - 3}
      \sum_{k = j + 1}^{n - 2} 1
    &= \frac{1}{3} n^3 - 3 n^2 + \frac{26}{3} n - 8
    \eqend{,}
  \nexteq
      \underbrace{
        \npcauset[Chain2ChainlowerDupperB]{1/wild,2/left label=$i$,8,3/wild,4/left label=$j$,9,5,6,7/wild,10/left label=$k$,11/wild},
        \npcauset[Chain2ChainlowerBupperU]{1/wild,2,6,9,3,4,5/wild,7,8/wild,10,11/wild}
      }:
    &&
      2 \sum_{i = 1}^{n - 6}
      \sum_{j = i + 1}^{n - 5}
      \sum_{k = j + 3}^{n - 2} 1
    &= \frac{1}{3} n^3 - 5 n^2 + \frac{74}{3} n - 40
    \eqend{.}
  \end{align}
  \end{subequations}
\item The elements are linked together and each is linked to the chain twice, 
  \begin{align}
  \label{eq:PosetNumbers.MaxBut2Layers.Linked2Links2Links}
      \npcauset[Chain2ChainlowerBupperB]{1/wild,2/left label=$i$,6,3/wild,4/left label=$j$,9,5/wild,7/left label=$k$,8/wild,10/left label=$l$,11/wild}:
    &&
      \sum_{i = 1}^{n - 5}
      \sum_{j = i + 1}^{n - 4}
      \sum_{k = j + 1}^{n - 3}
      \sum_{l = k + 1}^{n - 2} 1
    &= \frac{1}{24} n^4 - \frac{7}{12} n^3 + \frac{71}{24} n^2 - \frac{77}{12} n + 5
    \eqend{.}
  \end{align}
\end{enumerate}
The poset number $p(n, n - 2)$ as given in \eqref{eq:PosetNumbers.MaxBut2Layers} is the sum of all the counts from these cases. 

Among the $p(n, n - 2)$ posets with cardinality $n$, there are the following locally symmetric posets. 
First, any locally unsymmetric poset $P$ with cardinality $n - 1$ and $n - 2$ layers form $n - 1$ equivalence classes in $P / \npcauset[Singleton]{1}$, and an additional element can be added to any of these classes so that there are $(n - 1) u(n - 1, n - 2)$ combinations. 
Second, any singleton-symmetric poset $P$ that is retractable to the $(n - 2)$-chain forms $n - 2$ equivalence classes in $P / \npcauset[Singleton]{1}$, and the remaining two elements are distributed to the classes chosen with repetition, which are $\binom{2 + (n - 2) - 1}{2}$ possibilities. 
And third, all $n - 3$ posets of the type $C_1 \sqcup \left( \npcauset[2Chain_2Chain]{3,4,1,2} \right) \sqcup C_2$ (with chains $C_1, C_2$ such that $\lvert C_1 \rvert + \lvert C_2 \rvert = n - 4$) are 2-chain-symmetric and retractable to the $(n - 2)$-chain. 
Together, we have 
\begin{align}
\label{eq:PosetNumbers.MaxBut2Layers.LocallySymmetric}
    \frac{1}{2} (n - 3) (n - 2) (n - 1)
    + \frac{1}{2} (n - 2) (n - 1)
    + n - 3
  &= \frac{1}{2} \left( n^3 - 5 n^2 + 10 n - 10 \right)
\end{align}
locally symmetric posets with cardinality $n$ and $n - 2$ layers. 
The number $u(n, n - 2)$ in \eqref{eq:PosetNumbers.MaxBut2Layers} is the difference between the previously computed number $p(n, n - 2)$ and this number of locally symmetric posets. 
\qed

The ratio of locally unsymmetric posets among the number of posets with $n - 2$ layers 
\begin{align}
\label{eq:PosetNumbers.MaxBut2Layers.Ratio}
    \frac{u(n, n - 2)}{p(n, n - 2)}
  &= 1 - \frac{3}{n} + \frac{3}{n^2} + \Ord\bigl( n^{-3} \bigr)
\end{align}
shows that almost all posets with $n - 2$ layers are locally unsymmetric for large $n$. 
Though, this ratio is smaller than the ratio \eqref{eq:PosetNumbers.MaxBut1Layers.Ratio} for all $n > 3$.

\subsection{Posets with a small number of layers}
\label{subsec:Enumeration.FewLayers}

Similar calculations could be conducted to get the poset numbers $p(n, n - 3)$, $p(n, n - 4)$, and so on, together with the respective numbers for the locally unsymmetric posets.
More relevant applications, however, are the ratios for 2-layer, 3-layer, 4-layer (and higher layer) posets.
In the large $n$ behaviour, these posets dominate and it was shown that a certain type of 3-layer posets suffice to get the asymptotic enumeration of all posets \cite{KleitmanRothschild:1975}.
Such posets are commonly referred to as Kleitman-Rothschild orders especially in the context of causal set theory~\cite{LoomisCarlip:2018}.
An example is given in \cite[Fig.~1]{CarlipCarlipSurya:2023}, which has multiple local symmetries and is retractable to the 3-chain,
\begin{subequations}
\eqseqlabel{eq:KROrders.ExampleContractions}
\begin{align}
      \nrcauset[KRExample]{6,4,14,2,12,1,10,9,8,7,16,5,15,3,13,11}{1/8,1/9,2/5,2/7,2/10,4/7,4/10,4/12,5/13,6/8,6/9,7/13,7/15,8/11,8/16,9/11,9/16,10/13,10/15,12/15}
    \oslash \npcauset[2Diamond]{1,3,2,4}
    \oslash \npcauset[2Diamond]{1,3,2,4}
    \oslash \npcauset[Singleton]{1}
  &= \npcauset[3DiamondLinked2Diamond]{4,8,7,1,5,9,3,2,6}
    \oslash \npcauset[2Diamond]{1,3,2,4}
    \oslash \npcauset[Singleton]{1}
  = \npcauset[2Diamond]{1,3,2,4}
    \oslash \npcauset[Singleton]{1}
  = \npcauset[3Chain]{1,2,3}
  \eqend{,}
\nexteq
      \nrcauset[KRExample]{6,4,14,2,12,1,10,9,8,7,16,5,15,3,13,11}{1/8,1/9,2/5,2/7,2/10,4/7,4/10,4/12,5/13,6/8,6/9,7/13,7/15,8/11,8/16,9/11,9/16,10/13,10/15,12/15}
    \oslash \npcauset[Singleton]{1}
    \oslash \npcauset[2Diamond]{1,3,2,4}
  &= \nrcauset[3Simplex012Faces]{6,4,12,2,10,1,8,7,14,5,13,3,11,9}{1/8,2/5,2/7,4/7,4/10,5/11,6/8,7/11,7/13,8/9,8/14,10/13}
    \oslash \npcauset[2Diamond]{1,3,2,4}
  = \npcauset[3Chain]{1,2,3}
  \eqend{.}
\end{align}
\end{subequations}

Since the number of layers represents ``instances of time'' in the context of causal set theory, most posets with a small number of layers have to be excluded when we want to approximate spacetime manifolds with these discrete structures.
So suppression mechanism are considered to remove the dominating low layer posets in favour for posets that may resemble discretized spacetime manifolds \cite{LoomisCarlip:2018,CarlipCarlipSurya:2023,CarlipCarlipSurya:2023b}.

\begin{remark}
  Spacetime manifolds, as the continuous counterparts to causal sets, are partially ordered by their causal structure.
  They are not locally finite, instead the interval $J^+(a) \cap J^-(b)$ between an ordered pair of points $a < b$ is a compact subset.
  Also they do not have local symmetries as they are distinguishing, meaning any two distinct elements are uniquely identified by their pasts and futures~\cite{MinguzziSanchez:2008}.
  This is similar to locally unsymmetric causal sets, where all elements are uniquely identified by their relations to all other elements in the poset.
\end{remark}

Getting an exact quantitative statement about the number of posets with few layers is far more complicated than for the poset numbers computed in the sections above. 
For an estimate of 2-layer posets, we consider a generic 2-layer poset with cardinality $n = n_1 + n_2 \geq 2$ such that $n_1 \geq 1$ elements are on the first layer and $n_2 \geq 1$ elements are on the second layer. 
All second layer elements have to cover at least one element on the first layer so that the least number of links is $n_2$. 
In any 2-layer poset $P$ with exactly $n_2$ links, the second layer elements cover $k$ elements of the first layer where $1 \leq k \leq \min(n_1, n_2)$, so there are
\begin{align}
  &\quad
    1
    + \sum_{k = 2}^{\min(n_1, n_2)}
      \left(
        \sum_{i_1 = 1}^{%
          \lfloor n_2 / k \rfloor %
        }
        \sum_{i_2 = i_1}^{%
          \lfloor (n_2 - i_1) / (k - 1) \rfloor %
        }
        \dots
        \sum_{i_{k - 1} = i_{k - 2}}^{%
          \lfloor (n_2 - i_1 - \dots - i_{k - 2}) / 2 \rfloor %
        }
          1
      \right)
\nnexteq
\label{eq:PosetNumbers.2Layers.FewestLinks}
  &= 1
    + \underbrace{
      \left\lfloor \frac{n_2}{2} \right\rfloor
    }_{\mathclap{n_1, n_2 \geq 2}}
    + \underbrace{
      \sum_{i_1 = 1}^{%
        \lfloor n_2 / 3 \rfloor %
      }
        \left( \left\lfloor \frac{n_2 - i_1}{2} \right\rfloor - i_1 + 1 \right)
    }_{\mathclap{n_1, n_2 \geq 3}}
    + \cdots
\end{align}
possibilities, where $\lfloor \cdot \rfloor$ is the floor function. 
The first term of this sum (1 possibility) corresponds to the unique poset with $n_2$ elements all covering the same element on the first layer ($k = 1$). 
The second term ($k = 2$) enumerates the posets such that the second layer elements covering two elements on the first layer --- these are the partitions of $n_2$ into two subsets each with at least one element, so $\lfloor \frac{n_2}{2} \rfloor$ possibilities. 
This continues recursively up to $k = \min(n_1, n_2)$ partitions. 
If $n_1 \geq n_2$, the last term in \eqref{eq:PosetNumbers.2Layers.FewestLinks} is 1, for the unique poset with $n_2$ separated 2-chains parallel composed with an $(n_1 - n_2)$-antichain.
All of these posets are at least singleton-symmetric (or 2-chain-symmetric in the last case for $n_1 = n_2$) except for the two special cases \npcauset[2Chain]{1,2} and \npcauset[2Chain_Singleton]{3,1,2}. 

For any of these posets, adding one more link connects two separated components and may break some of the local symmetries, but may also introduce a new local symmetry. 
For example, 
\begin{align}
\label{eq:2LayerPosets.AddingLinks.Example}
    \left( \npcauset[Vee_2Chain]{4,5,1,3,2} \right)
  &\to \npcauset[W]{3,5,1,4,2}
\end{align}
breaks the singleton-symmetry but introduces a 2-chain-symmetry. 
A locally unsymmetric poset is obtained by adding links such that all local symmetries are broken, link by link without introducing new local symmetries. 
The quickest way to obtain a locally unsymmetric poset in that way is by starting from the poset that is a parallel composition of 2-chains (and a singleton) such that $n_2 = n_1$ (or $n_2 = n_1 - 1$) and adding links until we reach a poset that is a parallel composition $\bigsqcup \mathfrak{S}$ of a subset
\begin{align}
\label{eq:2LayerPosets.AddingLinks.LocallyUnsymmetric.Quickest}
    \mathfrak{S}
  &\subset \{\npcauset[Singleton]{1}\}
    \cup \left\{ F_{2 j} \bypred j \geq 1 \right\}
\end{align}
such that $\mathfrak{S}$ contains at least one of the fences $F_{2 j}$. 
For example, 
\begin{align}
\label{eq:2LayerPosets.AddingLinks.LocallyUnsymmetric.Example}
    \left( \npcauset[7_2Chains_Singleton]{15,13,14,11,12,9,10,7,8,5,6,3,4,1,2} \right)
  &\to \left( \npcauset[8FenceN2ChainSingleton]{15,13,14,10,12,9,11,6,8,4,7,2,5,1,3} \right)
  \eqend{.}
\end{align}
Any further link added may or may not reintroduce local symmetries.
If we continue adding links, we eventually arrive at a poset with a bipartite subset.
Recall that all bipartite posets~\eqref{eq:PosetClass.2Chain.1Chain} are singleton-symmetric.
When reaching a sufficiently large number of links, at least two second-layer elements will cover all first layer elements and any such poset is singleton-symmetric. 

While we have now shown that 2-layer posets with a sufficiently small or large number of links are almost always locally symmetric, this cannot be expected for the typical 2-layer posets considered in~\cite{KleitmanRothschild:1975}. 
The expression~\eqref{eq:PosetNumbers.2Layers.FewestLinks} can be approximated by an $n$-th order polynomial, but there are far more posets with two layers (including those with an arbitrarily large dimension). 
Consider posets with a cardinality of $4 n$, $2 n$ elements on each layer and each element linked to $n$ elements. 
For a singleton-symmetry in the second layer, the chance of an element to cover the same first layer elements as another second layer element is $2^{-n}$. 
Even if such posets may have $(Q, r)$-symmetries with posets $Q$ larger than the 2-chain, like 
\begin{subequations}
\begin{align}
    \nrcauset[4Gonal]{6,4,8,2,7,1,5,3}{2/5,4/7}
    \oslash' \npcauset[Wedge]{2,1,3}
  &= \npcauset[Vee]{1,3,2}
  \eqend{,}
\nexteq
    \nrcauset[2GonalCrown]{9,7,5,12,3,11,2,10,1,8,6,4}{2/8,3/6,3/10,5/8,5/11,7/10}
    \oslash' \npcauset[2Gonal]{2,1,4,3}
  &= \npcauset[2Gonal]{2,1,4,3}
  \eqend{,}
\end{align}
\end{subequations}
we have to expected that the numbers $u(n, 2)$ and $p(n, 2)$ are almost equal when $n$ becomes large. 
In general, we may assume that almost all posets with a fixed cardinality and fixed but small number of layers are locally unsymmetric, though, a conclusive statement on this goes beyond this work.

\section{Final remarks and applications of local symmetries}
\label{sec:Conclusion}

Starting from locally unsymmetric posets, we can classify all posets by extensions with local symmetries.
I finish with some remarks on local extensions for the construction of symmetry classes, followed by a discussion of local symmetries in causal sets.

\subsection{Local retractions and local extensions}
\label{subsec:LocalOperations}

With a retraction $P \oslash_r Q$ as constructed above, all $r$-cyclic, local symmetries described by the finite poset $Q$ are removed from a poset $P$.
For the characterisation of all posets, we may rather want to retract individual $(Q, r)$-symmetries from a $(Q, r)$-symmetry set of $A \isom Q$ in $P$, which we denote as $A^{(Q, r)}$.
Let $\Sigma$ be the subgroup of all automorphisms that leave the elements in $P \setminus A^{(Q, r)}$ fixed.
Taking a quotient of $P$ by the subgroup $\Sigma$ yields a retracted poset $\hat{P}$ that is isomorphic to $P$ up to a finite subset $\hat{A}$ as the retracted version of $A^{(Q, r)}$.
This \emph{local retraction} is written as $\hat{P} = (P, A^{(Q, r)}) \triangleright \hat{A}$, and in the Hasse diagram representation, we may use a different colour to denote the subset $A^{(Q, r)}$.
\begin{example}
  For example,
  \begin{align}
  \label{eq:Example.LocalRetraction}
      \nrcauset[3CrownBrokenLowHighSelLow]{6/wild,5/wild,2,8,1,7,4,3}{2/7}
      \mathbin{\triangleright} \npcauset[Singleton]{1}
    &= \nrcauset[3CrownBrokenHigh]{5,2,7,1,6,4,3}{2/6}
    \eqend{,}
  \end{align}
  is the local retraction of a two element subset $A^{(\smallSingleton, 2)} \isom \left( \npcauset[2Antichain]{2,1} \right)$ on the first layer to the singleton, $\hat{A} \isom \npcauset[Singleton]{1}$.
\end{example}

Note that such local retractions have an inverse operation by construction. 
We can \emph{locally extend} a poset $\hat{P}$ at a subset $\hat{A}$ by a $(\hat{A}, r)$-symmetry to a poset $P$ that is isomorphic to $\hat{P}$ up to a finite subset $A^{(Q, r)}$ with $\hat{A} \isom Q$ if $\hat{P} = (P, A^{(Q, r)}) \triangleright \hat{A}$ where $A^{(Q, r)}$ retracts to $\hat{A}$. 
For this inverse operation, we write $P = (\hat{P}, \hat{A}) \triangleleft A^{(Q, r)}$ and use again a different color or a label to indicate the corresponding subsets in the Hasse diagrams.
\begin{example}
  For the previous example, we have 
  \begin{align}
  \label{eq:Example.LocalExtension}
      \nrcauset[3CrownBrokenHighSelLow]{5/wild,2,7,1,6,4,3}{2/6}
      \mathbin{\triangleleft} \left( \npcauset[2Antichain]{2,1} \right)
    &= \nrcauset[3CrownBrokenLowHigh]{6,5,2,8,1,7,4,3}{2/7}
    \eqend{.}
  \end{align}
\end{example}

With these local retractions and local extensions, it becomes much more clear if and how a sequence of retractions by local symmetries is commutative. 

Any $(Q, r)$-retraction of a $(Q, r)$-retractable poset $P$ decomposes into local $(Q, r)$-retractions (that are prime).
If $Q = \npcauset[Singleton]{1}$, then the local retractions commute, because any equivalence class in $P / \npcauset[Singleton]{1}$ is only affected by a single local retraction.
\begin{example}
  Two local singleton-extensions of the crown at different elements can be locally retracted independently, 
  \begin{align}
  \label{eq:ExampleContractionDecomposition}
      \nrcauset[3CrownBrokenLowHigh]{6,5,2,8,1,7,4,3}{2/7}
      \oslash \npcauset[Singleton]{1}
    &= \left(
        \nrcauset[3CrownBrokenLowHighSelLowHigh]{6/wild,5/{wild,label=1},2,8,1,7,4/wild,3/{wild,label=2}}{2/7}
        \mathbin{\triangleright_1} \npcauset[Singleton]{1}
      \right)
      \mathbin{\triangleright_2} \npcauset[Singleton]{1}
    = \left(
        \nrcauset[3CrownBrokenLowHighSelLowHigh]{6/wild,5/{wild,label=1},2,8,1,7,4/wild,3/{wild,label=2}}{2/7}
        \mathbin{\triangleright_2} \npcauset[Singleton]{1}
      \right)
      \mathbin{\triangleright_1} \npcauset[Singleton]{1}
    = \nrcauset[3Crown]{4,2,6,1,5,3}{2/5}
    \eqend{.}
  \end{align}
\end{example}
Commutativity is no longer guaranteed for larger symmetries.
\begin{example}
  The order of retractions of the 2-chain-symmetries in the following poset yields different results
  \begin{subequations}
  \begin{align}
      \nrcauset[3CrownBroken3CrownBrokenDouble2ChainSel]{9/sel2,7/{sel2,label={\color{sel2}2}},4,12/sel2,3/sel2,11/sel2,2/sel1,1/{sel1,label={\color{sel1}1}},10,8/sel2,6/sel1,5/sel1}{1/6,1/8,2/5,2/8,3/5,3/6,4/10,7/11}
      \oslash \npcauset[2Chain]{1,2}
    &= \left( 
        \nrcauset[3CrownBroken3CrownBrokenDouble2ChainSel]{9/sel2,7/{sel2,label={\color{sel2}2}},4,12/sel2,3/sel2,11/sel2,2/sel1,1/{sel1,label={\color{sel1}1}},10,8/sel2,6/sel1,5/sel1}{1/6,1/8,2/5,2/8,3/5,3/6,4/10,7/11}
        \mathbin{\color{sel1}\triangleright_1} \npcauset[2Chain]{1,2}
      \right)
      \mathbin{\color{sel2}\triangleright_2} \npcauset[2Chain]{1,2}
    = \nrcauset[3Crown]{4,2,6,1,5,3}{2/5}
  \nexteq
    &\neq \left( 
        \nrcauset[3CrownBroken3CrownBrokenDouble2ChainSel]{9/sel2,7/{sel2,label={\color{sel2}2}},4,12/sel2,3/sel2,11/sel2,2/sel1,1/{sel1,label={\color{sel1}1}},10,8/sel2,6/sel1,5/sel1}{1/6,1/8,2/5,2/8,3/5,3/6,4/10,7/11}
        \mathbin{\color{sel2}\triangleright_2} \npcauset[2ChainSel]{1/sel1,2/sel1}
      \right)
      \mathbin{\color{sel1}\triangleright_1} \npcauset[2Chain]{1,2}
    = \npcauset[N]{2,4,1,3}
    \eqend{,}
  \end{align}
  \end{subequations}
  because the second local retraction on the second equation line also involves elements retracted by the first local retraction.
\end{example}

By using local extensions to break symmetries of a poset partially, and then retracting the partially broken symmetries followed by retractions of all other symmetries, we can construct larger locally unsymmetric posets.
\begin{example}
\label{eg:Fences.FromPolygonalPosets}
  Recall that the $n$-gonal poset is an $(\npcauset[2Chain]{1,2}, n)$-extension of the 2-chain. 
  Any fence of even cardinality $n$ can be constructed recursively from the $(n - 1)$-gonal poset, breaking the $(\npcauset[2Chain]{1,2}, n - 1)$-symmetry when locally extending by a 2-antichain and then retracting by the $(n - 2)$-fence and the singleton-symmetry, for example, 
  \begin{align}
      \nrcauset[5GonalSelLow]{8/wild,6,10,4,9,2,7,1,5,3}{2/5,4/7,6/9}
      \mathbin{\triangleleft} \left( \npcauset[2Antichain]{2,1} \right)
      \oslash \npcauset[N]{2,4,1,3}
      \oslash \npcauset[Singleton]{1}
    &= \npcauset[6Fence]{4,6,2,5,1,3}
    \eqend{.}
  \end{align}
\end{example}
It remains to prove if all connected, locally unsymmetric posets of a given cardinality $n$ (except for the $n$-chain poset) are constructable from smaller locally unsymmetric posets by a series of local extensions and retractions.

The notation for local extensions can be used also in the construction of symmetry classes for larger posets from symmetry classes of smaller ones.
\begin{example}
  All (prime) 2-chain extensions of the following 3-layer poset are
  \begin{align}
      \left[
        \npcauset[Nup]{2,5,1,3,4}
        \odot \npcauset[2Chain]{1,2}
      \right]'
    &\pickindent{=} \bigcup_{
        S \in [\smallTwoChain \odot \smallTwoChain]'
      }
      \left\{
        \npcauset[NupSel13]{2,5,1/wild,3/wild,4} \mathbin{\triangleleft} S,
        \npcauset[NupSel23]{2/wild,5,1,3/wild,4} \mathbin{\triangleleft} S,
        \npcauset[NupSel25]{2/wild,5/wild,1,3,4} \mathbin{\triangleleft} S,
        \npcauset[NupSel34]{2,5,1,3/wild,4/wild} \mathbin{\triangleleft} S
      \right\}
    \eqbreakr
      \cup \bigcup_{
        S_1, S_2 \in [\smallTwoChain \odot \smallTwoChain]'
      }
      \left\{
        \npcauset[NupSel2513]{2/{wild, label=2},5/wild,1/{wild, label=1},3/wild,4} \mathbin{\triangleleft_1} S_1 \mathbin{\triangleleft_2} S_2,
        \npcauset[NupSel2534]{2/{wild, label=2},5/wild,1,3/wild,4/{wild, label=1}} \mathbin{\triangleleft_1} S_1 \mathbin{\triangleleft_2} S_2
      \right\}
    \eqend{.}
  \end{align}
  Using~\eqref{eq:CountingFunction.2Chain.2Chain}, the enumeration of this class is
  \begin{align}
      c_n\left[
        \npcauset[Nup]{2,5,1,3,4}
        \odot \npcauset[2Chain]{1,2}
      \right]'
    &= \begin{cases}
        4 n - 28 & \text{if $n$ is odd and~} n \geq 11 \eqend{,} \\
        10 & \text{if~} n = 9 \eqend{,} \\
        4 & \text{if~} n = 7 \eqend{,} \\
        0 & \text{otherwise.}
      \end{cases}
  \end{align}
  To proof this via induction, note that the 2-chain extensions always add an even number of elements (so the count for even $n$ vanishes) and the first cases up to $n = 11$ can be checked directly.
  For the general induction step $\hat{n} = n + 2$, there are just 4 more posets with a double 2-chain extension (compared to the case $n$) for each of the two double extension possibilities.
  So there are $4 n - 28 + 8$ posets for $\hat{n} = n + 2$, which is $4 \hat{n} - 28$.
\end{example}
Similarly, the 2-chain extension class for any (finite) locally unsymmetric poset is constructed by replacing a single pair of linked elements, then two independent pairs of linked elements, then three, and so on.
For a 2-chain extension class of a locally symmetric poset, the equivalence classes under the respective local symmetries have to be respected so that elements from one equivalence class are only counted once per local 2-chain extension.

\subsection{Sprinkled causal sets and stability of locally unsymmetric posets}

To conclude this study of local symmetries in posets, consider causal sets that are generated from a given spacetime manifold via a Poisson process (\emph{sprinkling}) with a probability measure as follows.
For any (pre)-compact region $U$ of a spacetime manifold $M$ with volume measure $\nu$ and a sprinkling density $\rho > 0$ (a strictly positive real parameter), the probability to obtain a random sprinkle $\mathsf{S}_U$ (as a causal set) with cardinality $n$ is given by 
\begin{align}
\label{eq:SprinklingMeasure}
    \Pr\Bigl( \lvert \mathsf{S}_U \rvert = n \Bigr)
  &= \frac{\rho^n}{n!} \nu(U)^n \e^{-\rho \nu(U)}
  \eqend{.}
\end{align}
We can then use this measure to construct a probability space for an entire spacetime manifold, see \cite{FewsterHawkinsMinzRejzner:2021} for a mathematical review.
Especially for Minkowski spacetime, sprinkles are typically not only locally unsymmetric but exhibit an even stronger property defined as follows.

\begin{definition}
  For any $k \in \Naturals_0$, a poset $P$ is \defof{$k$-stable locally unsymmetric} if, for every subset $S \subseteq P$ that has cardinality $0 \leq \lvert S \rvert \leq k$, the poset $P \setminus S$ is locally unsymmetric. 
  A poset $P$ is \defof{total locally unsymmetric} if $P \setminus S$ is $k$-stable locally unsymmetric for every finite $k < \lvert P \rvert$. 
\end{definition}

\begin{example}
\label{eg:TotalLocallyUnsymmetric.Chains}
  Any chain poset is total locally unsymmetric, because chains are locally unsymmetric and any subset of a chain is again a chain.
\end{example}

\begin{table}
  \centering
  \begin{tabular}{r|*{7}r|r}
  \toprule
    $s(n, l)$
         &    1 &    2 &     3 &     4 &     5 &     6 &     7 & $s_n$ \\  
  \midrule
       1 &    1 &      &       &       &       &       &       &     1 \\  
       2 &    0 &    1 &       &       &       &       &       &     1 \\
       3 &    0 &    0 &     1 &       &       &       &       &     1 \\
       4 &    0 &    0 &     1 &     1 &       &       &       &     2 \\
       5 &    0 &    0 &     0 &     3 &     1 &       &       &     4 \\
       6 &    0 &    0 &     2 &     8 &     6 &     1 &       &    17 \\
       7 &    0 &    0 &     4 &    37 &    36 &    10 &     1 &    88 \\
  \bottomrule
  \end{tabular}
  \caption{\label{tab:PosetNumbers.StableLocallyUnsymmetric} Numbers of 1-stable locally unsymmetric posets by layers and in total.}
\end{table}
\begin{remark}
  For finite posets, the numbers of 1-stable locally unsymmetric posets $s(n, l)$ by layer are given in~\autoref{tab:PosetNumbers.StableLocallyUnsymmetric}. 
  In~\autoref{eg:TotalLocallyUnsymmetric.Chains}, we saw that all chains are total locally unsymmetric, so that also $s(n, n) = 1$.
  Furthermore, it is easy to see that (similar to the calculation in~\autoref{subsec:Enumeration.MaxBut1Layers}), $s(n, n - 1) = \frac{1}{2} (n - 2) (n - 3)$.
  Note that up to cardinality $n = 7$, there is only one 2-layer poset (the 2-chain) that is 1-stable locally unsymmetric.
\end{remark}

\begin{theorem}
  A sprinkle in $d$-dimensional Minkowski spacetime is total locally unsymmetric with probability 1. 
\end{theorem}
\begin{figure}
  \centering
  \includegraphics{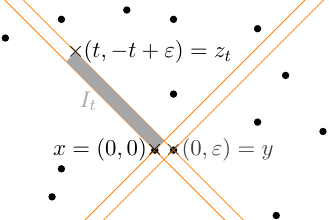}
  \caption{\label{fig:SprinkleTotalLocallyUnsymmetric} Choice of coordinates and construction of the region $I_t$ for a pair of separated elements in a sprinkle in Minkowski spacetime.}
\end{figure}
\proof
In $\mathbb{M}^{1 + 0}$, every sprinkle is a chain poset and hence total locally unsymmetric.
For any higher dimension $d > 0$, we prove it by contradiction and assume that a sprinkle $\mathsf{S}$ in $\mathbb{M}^{1 + d}$ contains a pair of distinct elements $x \neq y \in \mathsf{S}$ that has the same future up to a finite set of elements $F$, and the pair becomes a local symmetry in $\mathsf{S} \smallsetminus F$.
We choose a coordinate system such that $x$ is at the origin and $y$ has the same coordinates except for its first spatial coordinate that is displaced by a distance $\varepsilon > 0$.
For any $t > 0$, let $z_t$ be the point $x$ shifted by $t > 0$ in time and $- t + \varepsilon$ in the first spatial coordinate.
Let $I_t$ denote the spacetime interval from $x$ to $z_t$ as shown with an example for $\mathbb{M}^{1 + 1}$ in Fig.~\ref{fig:SprinkleTotalLocallyUnsymmetric}.
There are $n \geq 0$ sprinkled elements that fall in the compact region $I_t$ (which are also elements in $F$ by assumption) with probability
\begin{equation}
\label{eq:Sprinkling.ProbabilityMeasure.EmptyRegion}
    \Pr\Bigl( \lvert \mathsf{S} \cap I_t \rvert = n \Bigr)
  = \mathrm{e}^{-\rho \nu(I_t)} \frac{\rho^n}{n!} \nu(I_t)^n
  \eqend{.}
\end{equation}
The volume $\nu(I_t)$ does not vanish since $\varepsilon$ can be arbitrarily small but is fixed to a strictly positive value, and $\nu(I_t)$ increases with $t$ linearly for $d = 1$ and more than linearly for $d > 1$ by construction. 
Because we can choose $t$ arbitrarily large, consider $t \to \infty$, where the volume diverges and the probability~\eqref{eq:Sprinkling.ProbabilityMeasure.EmptyRegion} becomes zero for any finite $n$. 
This implies that the region $I_t$ has to contain an infinite number of elements with probability 1, which contradicts the initial assumption that $F$ is finite.
So we would have to remove an infinite number of elements to make $x$ and $y$ locally symmetric.
Hence any random sprinkle is total locally unsymmetric with probability 1.
\qed
\begin{corollary}
	A sprinkle in $d$-dimensional Minkowski spacetime does not contain singleton-symmetric elements (also known as \emph{non-Hegelian} pairs~\cite{WuethrichCallender:2015}) with probability 1.
\end{corollary}

Spacetime manifolds that are considered relevant in physics are \emph{globally hyperbolic} or at least \emph{distinguishing}, meaning that any two distinct points also have distinct pasts (and futures)~\cite{MinguzziSanchez:2008}.
Since any pair of elements of a locally unsymmetric (and locally finite) poset is distinct by their relations to all other elements, posets that are locally unsymmetric can be taken as a generalisation of distinguishing spacetimes.
The property of total local unsymmetry may suffice to rule out infinite posets with a finite number of layers, although this is not shown here.
Further investigations are necessary to determine if general (total) locally unsymmetric causal sets are a candidate for ``manifold-like'' causal sets and thus relevant for applications in quantum gravity, which provides an alternative to previous works~\cite{MajorRideoutSurya:2009}.

It is also an open question, whether local symmetries in finite posets may have direct applications in other parts of causal set theory, in particular, the (covariant) growth dynamics, see~\cite{Zalel:2023}.

\section*{Acknowledgement}
Thanks goes to Steven Carlip and Sumati Surya for discussions about Kleitman-Rothschild orders during the Quantum Gravity conference 2023 at the RRI, India, and to Yasaman Yazdi for suggestions on the presentation. 
Credits go to Eric Dolores for finding a typographical mistake in an earlier version.

\printbibliography

\end{document}